\documentstyle[11pt]{article}

\newtheorem{theorem}{{\sc Theorem}}
\newtheorem{proposition}{{\sc Proposition}}
\newtheorem{lemma}{{\sc Lemma}}
\newtheorem{corollary}{{\sc Corollary}}

\newtheorem{definition}{{\sc Definition}}

\newenvironment{proof}{\bigskip \noindent
         {\bf Proof.}}{ \hfill \fbox{\mbox{}} \medskip}

\newcommand{\cC}{{\cal C}}

\newcommand{\cL}{{\cal L}}

\newcommand{\cV}{{\cal V}}

\newcommand{\cZ}{{\cal Z}}

\newcommand{\bR}{{\bf R}}
\newcommand{\bZ}{{\bf Z}}

\newcommand{\bP}{{\bf P}}
\newcommand{\bE}{{\bf E}}

\newcommand{\bQ}{{\bf Q}}

\newcommand{\sip}{\beta} 

\newcommand{\philh}{J_n^{\cL_{1/2}}} 
 
\newcommand{\philp}{J_n^{\cL_{+}}}  
\newcommand{\philm}{J_n^{\cL_{-}}}  
\newcommand{\philpr}{J_n^{\cL_{1/2+r}}}

\begin{document}
 
\title{\bf Proof of the Conjecture that the Planar Self-Avoiding Walk
has Root Mean Square Displacement Exponent $3/4$}
\author{\sc Irene Hueter\footnotemark[1]} 
\date{} 

\maketitle

\renewcommand{\thefootnote}{\fnsymbol{footnote}} 
\footnotetext[1]{Department of Mathematics, University
                 of Florida, PO Box 118105, Gainesville, FL 32611-8105,
                 USA, email: {\sf hueter@math.ufl.edu.}} 
\footnotetext[2]{Date: July 25, 2001.}  
\footnotetext[3]{2000 Mathematics Subject Classification: 
              {\em Primary:} 60G50. 
              {\em Secondary:} 60D05, 60G57.}
\footnotetext[4]{Key words and phrases: 
         mean square displacement exponent,
         Palm distribution,  
         point process of self-intersections, 
         random polymer, 
         self-avoiding walk,
         self-intersection local time.}


\renewcommand{\thefootnote}{arabic{footnote}}

\vspace*{-0.8cm} 
\begin{abstract}
This paper proves the long-standing open conjecture
rooted in chemical physics ({\sc Flory} (1949)
\cite{flor}) that the {\em self-avoiding walk} (SAW)
in the square lattice has root mean square displacement 
exponent $3/4.$ We consider 
(a) the point process of self-intersections defined via 
certain paths of the symmetric simple random walk 
in ${\bf Z}^2$ and 
(b) a ``weakly self-avoiding cone process'' relative to this
point process when in a certain ``shape". 
We derive results on the asymptotic expected 
distance of the {\em weakly} SAW with 
parameter $\sip >0$ from its starting point, 
from which a number of distance exponents
are immediately collectable for the SAW as well. 
Our method employs the {\em Palm distribution} of the 
point process of self-intersection points in a cone.
\end{abstract} 

\pagestyle{myheadings}
\thispagestyle{plain}
\markboth{\sc Irene Hueter}
{\sc The Planar Self-Avoiding Walk}

\section{Introduction}

The self-avoiding walk serves as a model for linear polymer molecules.
Polymers are of interest to chemists and physicists and are
the fundamental building blocks in biological systems. A polymer
is a long chain of monomers (groups of atoms) joined to one another
by chemical bonds. These polymer molecules form together randomly
with the restriction of no overlap. This repelling force makes
the polymers more diffusive than a simple random walk.
What are the properties of an average configuration of a polymer,
most of all, what is the average distance between the two ends 
of a long polymer ? This paper will give an answer to the latter
question in two dimensions and a separate paper will deal with
all dimensions. Note that the subject ``self-avoiding walk" 
has been well attended, not exclusively by probabilists. Here,
no attempt is made to survey the vast
literature (consult e.g.\ {\sc Madras and Slade} \cite{masl}).
 
\smallskip
{\bf (Weakly) Self-Avoiding Walk.}
We will consider the weakly self-avoiding walk in ${\bf Z}^2$ starting 
at the origin. More precisely, if $J_n= J_n(\cdot)$ denotes the
{\em number} of {\em self-intersections} of a symmetric simple random walk 
$S_0={\bf 0}, S_1, \ldots, S_n$ on the planar lattice 
starting at the origin, that is,
\begin{equation}
  \label{selfintersections}
      J_n = J_n(S_0, S_1, \ldots, S_n) = \sum_{0 \leq i < j \leq n} 
            1_{\{ S_i = S_j \}} ,
\end{equation} 
and if $\sip \geq 0$ denotes the self-intersection parameter,
then the {\em weakly self-avoiding walk} is the stochastic 
process, induced by the probability measure 
\begin{equation}
  \label{intersectmeasure}
    \bQ^{\sip}_n(\cdot) = \frac{\exp \{ - \sip J_n(\cdot) \}}{
                              \bE \exp \{ - \sip J_n(\cdot) \}} \, ,
\end{equation}
where $ \bE$ denotes expectation of the random walk.
In other words, 
$J_n = r$ self-intersections are penalized by the
factor $\exp \{ - \sip r \}.$ The measure $ \bQ^{\sip}_n$ 
may be looked at as a measure on the set of all 
simple random walks of length $n$ which weighs relative to
the number of self-intersections. This restraint walk is 
also being called the {\em Domb-Joyce model} in the 
literature (see {\sc Lawler} \cite{lawl}, p.\ 170),
but importantly, differs from the discrete {\em Edwards model}, 
which is a related repelling walk (see {\sc Madras
and Slade} \cite{masl}, 
p.\ 367 and {\sc Lawler} \cite{lawl}, p.\ 172 for some background). 
Thus, when setting $\sip=0,$ we recover the simple random 
walk (SRW), while letting $\sip \rightarrow \infty$ well mimics
the {\em self-avoiding walk} (SAW). The SAW 
in ${\bf Z}^2$ is a SRW-path of length $n$
with {\em no} self-intersections. 
Thus, this walk visits each site of its path {\em exactly} once.

We will be interested in the expected distance of the weakly
SAW from its starting point after $n$ steps,
as measured by the root mean square displacement
at the $n$-th step. Let $\bE_{\sip} = \bE_{\bQ^{\sip}_n}$
denote expectation under the measure $\bQ^{\sip}_n$
(that is, expectation wrt.\ to the weakly SAW). Thus, 
$\bE_0$ denotes expectation wrt.\ to the SRW. Also, write  
$S_n = (X_n, Y_n)$ for every integer $n \geq 0.$ 
Objects of interest to us are the expectation $\bE_{\beta}$
of the {\em distance} 
$$
  \chi_n = (X_n^2 + Y_n^2)^{1/2} 
$$ 
of the walk from the starting point ${\bf 0},$  
the {\em mean square displacement} 
$\bE_{\beta} \chi_n^2,$ and the 
{\em root mean square displacement} 
$ ( \bE_{\beta} \chi_n^2 )^{1/2}$
of the weakly SAW. Shorter, we shall write MSD and RMSD (for the
latter two), respectively. The main results of this paper are
the following statements, valid in two dimensions.

\begin{theorem}
 \label{sawdistance}
For each $\sip >0,$  
the exponent of the distance of the planar
weakly self-avoiding walk equals $3/4.$  
Furthermore, there are some constants $0 < \rho_1 = \rho_1(\sip) \leq  
 \rho_2 < \infty$ ($\rho_2:$ uniform in $\sip$) such that 
$$  
   \rho_1  \leq
     \liminf_{n \rightarrow \infty}
       n^{- 3/4} \,  \bE_{\sip} (\chi_n) 
        \leq  \limsup_{n \rightarrow \infty}
       n^{- 3/4} \,  \bE_{\sip} (\chi_n) 
         \leq \rho_2.
$$ 
In particular,
the planar self-avoiding walk has distance exponent $3/4.$ 
\end{theorem}

The proof of Theorem \ref{sawdistance} 
is collected in Propositions \ref{msdupper} and \ref{msdlower}
and in Corollary \ref{liminfsupmsd}. We remark that
the constants $\rho_2$ and $\rho_4$ in Theorems
\ref{sawdistance} and \ref{rmsd}, respectively, are independent
of all $\sip >0,$ whereas the constants $\rho_1$ and $\rho_3,$ 
respectively, in the lower bounds may depend on $\sip,$ even as 
$\sip \rightarrow 0$ or $\sip \rightarrow \infty.$

\begin{theorem}
 \label{rmsd}
For each $\sip >0,$ 
the root mean square displacement exponent of the planar
weakly self-avoiding walk equals $3/4.$ 
Moreover, there are some constants 
$0 < \rho_3= \rho_3(\sip) \leq \rho_4 < \infty$ 
($\rho_4:$ uniform in $\sip$) such that 
$$
   \rho_3  \leq
     \liminf_{n \rightarrow \infty}
       n^{- 3/2} \, \bE_{\sip} (\chi_n^2) 
        \leq  \limsup_{n \rightarrow \infty}
       n^{- 3/2} \, \bE_{\sip} (\chi_n^2) 
         \leq \rho_4 .
$$ 
Particularly,
the planar self-avoiding walk has root mean square displacement
exponent $3/4.$ 
\end{theorem}

See Corollary \ref{variance} for a proof.
Theorem \ref{rmsd} solves a major several decades-old 
open conjecture that can be traced back to at least
{\sc Flory}'s work \cite{flor} in the 1940ies and
is one among numerous stones yet to be uncovered in the 
field of $2,3,$ and $4$-dimensional random polymers.

In contrast to popular believes, the same approach, 
the approach that we pursue 
in this article extends to dimensions $3, 4,$ and 
higher as well as $1.$ These cases
will be discussed in detail in a separate paper \cite{huet1} to come. 
Here, we content ourselves with stating the formula for
the RMSD exponent $\nu$ of the weakly SAW in ${\bf Z}^d$
(with all definitions being much the same as in two dimensions).
This expression coincides with the one for the RMSD exponent 
of the SAW, defined
as the $\lim_{n \rightarrow \infty} \, \lim_{\sip \rightarrow \infty}
       \, \ln  (\bE_{\beta} ( \chi_n^2))/ (2  \ln n)$
if the limits exist, in view of existing
thresholds of $n$ that are uniform in $\sip$ as $\sip \rightarrow
\infty$ and considerations towards exchanging the limits
$\lim_{\sip \rightarrow \infty}$ and
$\lim_{n \rightarrow \infty}:$
\begin{eqnarray*}
        \nu & = & 1  \qquad \qquad \qquad \qquad \mbox{ for } d =1, \\
             & = &   \max( \frac{1}{2}, \frac{1}{4} + \frac{1}{d})
                  \qquad \mbox{ for } d \geq 2. 
\end{eqnarray*}
Equivalently, 
\begin{eqnarray*}
        \nu & = &  1  \quad \, \, \, \, \, \; \, \mbox{ for } d =1, \\
            & = &  3/4  \, \, \quad \mbox{ for } d =2, \\
            & = &  7/12  \quad \mbox{ for } d =3, \\
            & = &  1/2  \, \, \, \quad \mbox{ for } d \geq 4.
\end{eqnarray*}

We observe that, in dimension $4,$ the exponent $1/2$ arises
for different reasons than it occurs in dimensions $5$ and
higher. In dimension larger than $1,$ 
the MSD is asymptotic to the sum of at least two terms, one of which
is dominating in dimensions $2$ and $3,$ the other
of which is dominating for $d \geq 5.$ The latter is the
term that would present itself for the SRW. 
In this sense, the SAW in dimensions $d \geq 5$ behaves as the SRW.
Note that $7/12=0.58333...$ differs from the value $0.59...$ (see e.g.\
{\sc Lawler} \cite{lawl}, p.\ 167) that was believed more recently,
as stipulated by heuristic and ``numerical evidence" 
(Earlier estimates included the Flory estimate $0.6$).

Whereas our result is novel for $d=2,3,$ and $4$ and 
$\sip \in (0, \infty],$
the result on the RMSD exponent for the SAW for $d \geq 5$ 
is in {\sc Hara and Slade} \cite{hasl}
and the one on the RMSD exponent for the weakly SAW for $d=1$ 
is in {\sc Greven and den Hollander} \cite{grho}.
The former was accomplished via the perturbation technique
``lace expansion" and the latter via large deviation theory.
{\sc Brydges and Spencer} \cite{brsp} establish that the
scaling limit of the weakly SAW is Gaussian for sufficiently 
small $\sip >0 $ and $d \geq 5.$

Here are a few words about the ideas of proof and how we came
across them before we embark on the detailed arguments.
The square root of the mean square displacement of the SRW up
to time $n$ is of order $\sqrt{n}.$ Similarly, the scaling 
$\sqrt{t}$ is an elementary and distinctive feature of standard 
Brownian motion in ${\bf R}^d$ run for time $t.$ 
The latter may be calculated
by integrating the appropriate expression in polar coordinates,
in other words, by regarding the projection of the Brownian
motion onto any fixed line. This process along a line
is {\em one}-dimensional. Motivated by
this observation, it is natural to similarly attempt to
view a one-dimensional process that mimics the SRW and
which is penalized according to the number of self-intersections 
of the random walk which happen on (or near) a line 
and to ask the questions
(i) how far such a process is expected to move from the
origin and 
\mbox{(ii) how} its expected distance compares to
the one of the weakly SAW.
Clearly, if we fix a line, the random walk that we run
may not intersect it -- a configuration that is far from ideal --.
The situation improves if we pick a ``typical'' line. 
Of course, such a typical line would have to be 
chosen differently for each realization of a SRW-path. 
A useful concept in stochastic geometry 
allows us to deal with typical random geometric objects,
for instance, sample points of point processes. 
This is the {\em Palm distribution} of a random measure.

We shall employ the Palm distribution
of the point process of self-intersections, defined
via certain paths of length $n$ of the symmetric
SRW in ${\bf Z}^2,$ in a cone to introduce 
a ``weakly self-avoiding cone process'' relative 
to this point process when in a certain ``shape''.
To finish the story all the way,
the asymptotic expected distance of this
process can be calculated rather explicitly as n 
tends large, at least if the point process is 1/2-shaped, 
in which case it can be shown to equal the expected 
distance of the weakly SAW from its starting point. 
From these results along with some considerations
towards uniform bounds and estimates in $\sip$ as 
$\sip \rightarrow \infty,$ the mean square displacement 
exponent of the SAW immediately derives.

Section 2 presents a characterization of the SRW-paths that are
``atypical" but significant from the perspective of a weakly SAW.
Section 3 makes a connection between Palm distributions 
and the random walk, introduces this weakly self-avoiding cone 
process, calculates some asymptotic mean distances of this 
process and links those to the ones of the weakly SAW. 
Some remarks on the transitions
$\sip \rightarrow \infty$ and $\sip \rightarrow 0$ end Section 3.


\section{Atypical SRW-Paths are the Important Ones}
\setcounter{equation}{0} 

An elementary calculation, based on the Local Central Limit 
theorem (consult any graduate-level probability textbook),
shows that the expected {\em self-intersection local time}
(SILT) $\bE_0 J_n$ of the SRW is asymptotic to
$ \pi^{-1} n \ln n, $ with the error being no larger than order $n.$
However, these typical paths are ``negligible" in any analysis of
the weakly SAW as the following estimates indicate. 
We will make use of the convenient $o(\cdot)$ notation, that is,
write $f(n)= o(g(n))$ as $n \rightarrow \infty$ for two real-valued
functions $f$ and $g$ if $\lim_{n \rightarrow \infty} f(n)/g(n) =0.$
 
\begin{proposition}
   \label{upperboundforjn}{\bf (Upper Bound for $J_n$)}
Let $\beta > 0$ and let $\nu_0$ denote
the exponent of the number of self-avoiding walks.
Then for every $B > B_* = (\ln 4 - \nu_0)/\sip >0$ 
and every integer $n \geq 0,$
$$
   \bE_0 ( e^{ - \sip J_n} \, 1_{\{ J_n > B n \}} ) < 
            \bE_0 ( e^{ - \sip J_n} \, 1_{\{ J_n =0 \}} ),  
$$
in particular, as $n \rightarrow \infty,$
$$
   \bE_0 ( e^{ - \sip J_n} \, 1_{\{ J_n > B n \}} ) = 
            o( \bE_0 ( e^{ - \sip J_n} \, 1_{\{ J_n =0 \}} )) .  
$$
\end{proposition}

\begin{proof}
We proceed to show that the set of self-avoiding paths of the SRW
contributes a term to $ \bE_0 \exp \{ - \sip J_n  \}$ which
has larger exponential rate than the one contributed by the SRW
paths with $J_n  > B n $ for every $B >B_*.$

For this purpose, let $\Gamma_n$ denote the set of SAW-paths
$S_0={\bf 0}, S_1, \ldots, S_n$ (with $J_n =0$) up through time $n.$
Since, for every pair $(m,n),$ concatenating two paths 
$\gamma_1 \in \Gamma_n$ and $\gamma_2 \in \Gamma_m$ does not 
always provide a path in $\Gamma_{n+m},$ we immediately gain
$$
  \vert \Gamma_{n+m} \vert \leq \vert \Gamma_n \vert \, 
                                \vert \Gamma_m \vert,
$$
where $\vert \Gamma_n \vert $ denotes the cardinality of 
$\Gamma_n.$ An easy subadditivity argument yields that the
limit 
\begin{equation}
  \label{subaddlimit}
   \lim_{n \rightarrow \infty} \vert \Gamma_n \vert^{1/n}
   = e^{\nu_0}
\end{equation}
exists for some $2 \leq  e^{\nu_0} \leq 3$ and that
$ \vert \Gamma_n \vert \geq   e^{n \nu_0}$ 
for every integer $n \geq 0.$ Observe that in ${\bf Z}^d$
for $d >1,$ the {\em upper} bound $2d-1$ for $e^{\nu_0}$ may be
seen by counting all paths of length $n$ that do not return
to the  most recently visited point (clearly, an overestimate),
whereas the {\em lower} bound $d$ for $e^{\nu_0}$ may be seen
by counting all paths of length $n$ that take only positive
steps in both coordinates, e.g.\ for $d=2,$ 
move only north or east, say.

Fix $\sip >0$ and choose $B > (\ln 4 - \nu_0)/\sip >0.$
First, if $\bP_0$ denotes probability wrt.\ to the SRW,
\begin{eqnarray}
  \label{uppthresh}
  \bE_0 ( e^{ - \sip J_n} \, 1_{\{ J_n > B n \}} )  
         & < &  e^{ - \sip B n} \, \bP_0 ( J_n > B n)
              \leq  e^{ - \sip B n} .
\end{eqnarray}
Second, in view of the illustrated submultiplicativity 
property of the SAW,
\begin{eqnarray}
   \label{lowthresh}
    \bE_0 ( e^{ - \sip J_n} \, 1_{\{ J_n =0 \}} ) 
          & = & \bP_0 ( J_n = 0)   \\*[0.15cm] 
          & \geq &  e^{n \nu_0} 4^{-n} 
                   = e^{ -n (\ln 4 - \nu_0)}.   \nonumber
\end{eqnarray} 
Combining (\ref{uppthresh}) and (\ref{lowthresh}) and recalling 
that $ - \sip B < -(\ln 4 -  \nu_0 )$ yields
\begin{eqnarray*}
 \bE_0 ( e^{ - \sip J_n} \, 1_{\{ J_n > B n \}} )   
          & < & e^{ - \sip B n}  \\*[0.15cm]
          & < & e^{ -n (\ln 4 - \nu_0)} \\*[0.15cm] 
          & \leq  &  
              \bE_0 ( e^{ - \sip J_n} \, 1_{\{ J_n =0 \}} ), 
\end{eqnarray*}
and thus, both advertized claims. This ends the proof.
\end{proof}

Therefore, we learn that it suffices to focus on the 
SRW-paths that exhibit $J_n \leq n B_* = n (\ln 4 - \nu_0)/\sip.$
On the other hand, the paths with $J_n $ of order
less than $n$ are not significant either.

\begin{proposition}{\bf (Lower Bound for $J_n$)} 
    \label{lowerboundforjn} 
Let $\beta > 0.$ There is some $b_*= b_*(\sip) >0$ 
(made precise below)
so that for every $\delta >0$ and every $b < b_*,$
as $n \rightarrow \infty,$ 
$$
    \bE_0 ( e^{ - \sip J_n} \, 1_{\{ J_n \leq n^{1- \delta} \}} ) =
      o(   \bE_0 ( e^{ - \sip J_n} \, 1_{\{ J_n < b n \}} )).
$$       
\end{proposition}

\begin{proof}
As in the previous proof, it suffices to find an upper bound
for the lefthand side and a lower bound for the
righthand side of the display in such a fashion that 
the former is (exponentially) smaller than the latter. 

Fix $\delta>0.$ Let $\Omega^{\delta}_n$ denote
the set of all SRW-paths of length $n$ 
that have $J_n \leq n^{1 - \delta}.$
Clearly,
\begin{eqnarray}
     \bE_0 ( e^{ - \sip J_n} \, 1_{\{ J_n \leq n^{1- \delta} \}} ) 
            & \leq &  
             \bP_0 (J_n  \leq n^{1- \delta} )  \nonumber \\*[0.1cm]
          & = &  
     \label{thresholdup} 
            \exp \{ n ( \ln \vert \Omega^{\delta}_n \vert /n - \ln 4) \}.
\end{eqnarray} 

It remains to come up with a {\em lower} bound for 
$ \bE_0 ( e^{ - \sip J_n} \, 1_{\{ J_n < b n \}} )$ 
for all $b >0$ and to determine when this lower bound
is larger than the righthand side of (\ref{thresholdup}).
It will turn out that this is the case for all sufficiently
small $b.$ Note that if $b$ is not small enough, 
$e^{ - \sip J_n}$ may get too small. We begin with deriving
a lower bound for $\bP_0( J_n < bn ).$ 

Pick a suitable $s < b $ and 
consider the set $\Lambda_n$ of 
SRW-paths of length $n$ with $J_n = s n.$  
Observe that, for all sufficiently large $n,$
the set $\Lambda_n$ contains $\Omega^{\delta}_n.$
In fact, each path $\gamma \in \Omega^{\delta}_n$
gives rise to a large set $G_{\gamma}$ of paths in $\Lambda_n.$

To see this, choose any path $\gamma$ in $\Omega^{\delta}_n$
and introduce $\alpha_n n $
{\em repetitions} on that path, that is, choose $\alpha_n n $ 
distinct sites $x_j$ among the $n$ visited sites of $\gamma$ 
where $S_{j+1} =  S_{j-1}$ and 
$S_{j+2} =  S_j= x_j$ (Immediate backtracking and moving on).
If we choose $\alpha_n$ suitably, then this new SRW-path 
$\tilde{\gamma}$ in $G_{\gamma},$ starting at ${\bf 0},$
will have $J_n =sn.$ It will turn out that $\alpha_n = 
\tilde{\alpha} - \rho_n$ for some positive finite
constant $\tilde{\alpha}$ and $\rho_n \leq n^{- \delta}.$ 
Thus, $\alpha_n \rightarrow \tilde{\alpha}$
as $n \rightarrow \infty.$ To facilitate notation, we will
drop the subscript $n$ from $\alpha_n$ and just write $\alpha$
for the rest of the proof.
Observe that the path $\tilde{\gamma}$ of length $n$ 
has a trace which is shorter by $2 \alpha n$ units and ends
at the $n(1- 2 \alpha)$-th site of the path $\gamma,$ as 
$2\alpha n$ times were wasted revisiting sites. 
In order that two distinct paths $\gamma$ and $\gamma'$
in $\Omega^{\delta}_n$ generate two sets $G_{\gamma}$ and  
$G_{\gamma'}$ that are disjoint, 
first,
we do not allow to place repetitions among the $2 \alpha n$ last
sites of the paths, thus, let $\gamma$ in 
$\Omega^{\delta}_{n(1 - 2 \alpha)},$ 
and second, add a fixed number $f$ of self-intersection points 
at each site where there is at least one self-intersection point 
of $\gamma$ and accordingly deduct the corresponding
number of SILT from the $\alpha n $ repetitions to be performed.
Thus, in this repetition scheme, only one repetition is
allowed per site except for the prescribed repetitions
at existing self-intersection points of $\gamma,$
where a larger number of repetitions will be placed.
The latter agreement (concerning prescribed repetitions)
guarantees that the newly generated paths in
$G_{\gamma}$ distinguish themselves from all paths in $G_{\gamma'}$
for $\gamma$ different from $\gamma'$ and that
$G_{\gamma} \not = G_{\gamma'}$ (since each class can be uniquely
identified). In other words, the classes $G_{\gamma}$ are 
{\em disjoint.} Importantly, observe that the number of prescribed
repetitions is no larger than $f n^{1 - \delta} = o( sn)$ as 
$n \rightarrow \infty.$ Also, note that  
$\tilde{\gamma} \in  G_{\gamma}$ leave the same {\em trace}
as $\gamma  \in \Omega^{\delta}_{n(1 - 2 \alpha)}.$

Counting all paths in $ G_{\gamma}$ will lead to a lower
bound for $\vert \Lambda_n \vert,$ thus, 
to a lower bound for $\bP_0 (J_n = sn),$ 
and eventually, to a lower bound for $\bP_0 (J_n  \leq b n).$ 
Two moments' thoughts reveal
that the set of all selections of $\alpha n $
distinct sites among the $n(1 - 2 \alpha)(1+o(1))$ visited sites of 
$\gamma$ is in one-to-one correspondence with the set $ G_{\gamma}.$
The cardinality of the former equals the number of ways to 
distribute $\alpha n $ identical balls in $n(1 - 2 \alpha)(1+o(1))$ urns 
under the restriction that there is no more than one ball per urn.
(The correction $o(1)$ enters in view of the prescribed repetitions
at self-intersection points, which are negligible in number in comparison to 
the total $\alpha n,$ but will be suppressed in the calculation
below.) Hence, as an appeal to the Stirling approximation 
$k! = \sqrt{2 \pi k} e^{-k} k^k(1+o(1))$ as $k \rightarrow \infty,$  
\begin{eqnarray}
  \vert  G_{\gamma} \vert & = &
                   { n(1- 2 \alpha) \choose \alpha n } 
                   \nonumber  \\*[0.2cm]
                   & = & 
                   \label{cardggammaone} 
                  \left [  \frac{(1 - 2 \alpha)^{1 - 2 \alpha}}  
                    { \alpha^{\alpha} \, (1 - 3 \alpha)^{1 - 3 \alpha}} 
                    \right ]^n \, (\frac{1}{ 2 \pi n})^{1/2} \,
                    \left [ \frac{ 1 - 2 \alpha }{ 
                           \alpha  (1 - 3 \alpha)} \right ]^{1/2} 
                           \, (1 + o(1)) . 
\end{eqnarray}   
Since $(1 - 2 \alpha)^{1 - 2 \alpha}/(1 - 3 \alpha)^{1 - 3 \alpha} \geq 1,$
there is some $\xi = \xi_{\alpha} >0$ with 
$\exp \{ \xi \} \geq 1/\alpha^{\alpha}$ so that 
the righthand side of (\ref{cardggammaone}) 
$\geq \exp \{ \xi n \},$ that is, 
\begin{eqnarray}
    \label{cardggammatwo}  
      \vert  G_{\gamma} \vert & \geq & \exp \{ \xi n \} .
\end{eqnarray} 
We conclude that for each $\gamma \in 
\Omega^{\delta}_{n(1 - 2 \alpha)},$   
we can identify at least  $ \exp \{ \xi n \} $ SRW-paths with
$J_n = sn.$ Keeping in mind our earlier observations
that, first, all paths in $G_{\gamma}$ are different SRW-paths,
and second, the sets $G_{\gamma}$ are different from one another
and combining (\ref{cardggammatwo}) with the fact that each path in 
$\Omega^{\delta}_{n(1 - 2 \alpha)}$ can be completed to render
a path of length $n$ in $\Omega^{\delta}_n,$ we obtain
\begin{eqnarray}
  \vert \Lambda_n \vert & \geq & 
           \vert \Omega^{\delta}_{n(1 - 2 \alpha)} \vert
           \cdot \vert G_{\gamma} \vert
                 \nonumber  \\*[0.15cm]
           & \geq & 
              \label{largerthansawexpone} 
            \exp \{ n [ \ln \vert \Omega^{\delta}_{n(1 - 2 \alpha)} \vert/n
                    + \xi ] \} . 
\end{eqnarray}        
As a consequence, for each $b >0$ and suitably small $ 0<s< b,$ 
\begin{eqnarray}
   \bE_0 ( e^{ - \sip J_n} \, 1_{\{ J_n < b n \}} ) 
             & > & 
               e^{- \sip b n} \, \bP_0 ( J_n < bn ) 
              \nonumber  \\[0.15cm]
           & \geq &  
                e^{- \sip b n} \,  \bP_0 ( J_n = s n )  
             \nonumber \\[0.15cm]
             & = &  
                e^{- \sip b n} \, \vert \Lambda_n \vert \, 4^{-n}
             \nonumber \\*[0.15cm]
              & \geq &  
                 \label{jnlessthanbh} 
               \exp \{ n (
           \ln \vert \Omega^{\delta}_{n(1 - 2 \alpha)} \vert/n + \xi
            - \ln 4 - \sip b) \} .
\end{eqnarray} 
Thus, we arrived at a lower bound for 
$  \bE_0 ( e^{ - \sip J_n} \, 1_{\{ J_n < b n \}} ) .$

Next, we will verify the claim that
for $C =  2 \ln 4 +1$ and all sufficiently large $n,$
$$
  \frac{1}{n} \,  ( \ln \vert \Omega^{\delta}_{n} \vert -
   \ln \vert \Omega^{\delta}_{n(1 - 2 \alpha)} \vert)
     \leq C \alpha
$$
for every $\alpha.$ For this purpose, fix $\alpha.$ 
Each path in $\Omega^{\delta}_n$ arises either from a path in 
 $\Omega^{\delta}_{n(1 - 2 \alpha)} $ or from a path
in the difference $\Omega^{\delta'}_{n(1 - 2 \alpha)} \setminus
\Omega^{\delta}_{n(1 - 2 \alpha)}$ for $\delta  > \delta' \geq
\delta_*,$ where $\delta_*= \delta_*(n, \alpha)$ 
satisfies the equation $n^{1-\delta} = [n(1 - 2 \alpha)]^{1 - \delta_*}.$ 
Solving for $\delta_*$ gives
\begin{equation}
   \label{deltastar}
  \delta_* = \delta \, \frac{\ln n + \ln(1 - 2 \alpha)/\delta}{
              \ln n + \ln(1- 2 \alpha)}. 
\end{equation} 
Concatenating two paths 
$\gamma \in \Omega^{\delta}_{n(1 - 2 \alpha)}$ or $\in 
\Omega^{\delta'}_{n(1 - 2 \alpha)} \setminus 
\Omega^{\delta}_{n(1 - 2 \alpha)}$ and $\gamma' \in  
\Omega^{\delta}_{ 2 \alpha n}$ will not always produce
a path in  $\Omega^{\delta}_n$ since  $x^{1- \delta}$
is a concave function in $x>0,$ in particular, 
$(x_1 + x_2)^{1 - \delta} \leq x_1^{1 - \delta} +  x_2^{1 - \delta}$
for $x_1 , x_2 \geq 0,$
and, in the concatenated path, possible overlap of the subpaths 
$\gamma$ and $\gamma'$ may introduce additional SILT.
Consequently, we collect
\begin{eqnarray*}
  \vert \Omega^{\delta}_n \vert & \leq &
      \vert \Omega^{\delta}_{n(1 - 2 \alpha)} \vert \cdot
         \vert \Omega^{\delta}_{ 2 \alpha n} \vert  \, + \,
    \vert \Omega^{\delta'}_{n(1 - 2 \alpha)} \setminus
             \Omega^{\delta}_{n(1 - 2 \alpha)} \vert \cdot
         \vert \Omega^{\delta}_{ 2 \alpha n} \vert
                  \nonumber  \\*[0.15cm]
          & = &
    \vert \Omega^{\delta'}_{n(1 - 2 \alpha)} \vert \cdot
         \vert \Omega^{\delta}_{ 2 \alpha n} \vert \nonumber  \\
   \mbox{ equivalently, } 
      \qquad \qquad \qquad & &  \nonumber \\
   \ln \vert \Omega^{\delta}_n \vert 
      \, - \,  \ln \vert \Omega^{\delta'}_{n(1 - 2 \alpha)} \vert 
        & \leq &
       \ln \vert \Omega^{\delta}_{ 2 \alpha n} \vert
           \nonumber  \\*[0.1cm]
         & \leq &  2 \alpha n \ln 4.
\end{eqnarray*} 
Inspecting (\ref{deltastar}) makes clear that   
letting $n$ be suitably large will bring
$\delta_*$ as close to $\delta$ as desired.
Moreover, since the function  
$\ln \vert \Omega^{\delta}_n \vert/n$ is bounded and
nonincreasing in $\delta,$
it has at most finitely many jump discontinuities
of size larger than, say, $\tau >0.$ Therefore, in view of 
(\ref{thresholdup}) and the fact that $\vert \Omega^{\delta}_n \vert
\leq \vert \Omega^{\tilde{\delta}}_n \vert$ 
for $\tilde{\delta} < \delta,$ 
it suffices to restrict our attention to those $\delta$ for which
the function $\ln \vert \Omega^{\delta}_{n(1 - 2 \alpha)} \vert/n$ 
has no jumps of size  $\geq \alpha$ on the interval 
$[ \delta_*, \delta]$ for all sufficiently large $n.$
A combination of these observations leads us to conclude
that, for all sufficiently large $n,$ we obtain
$\ln (\vert \Omega^{\delta'}_{n(1 - 2 \alpha)}
     \vert /\vert \Omega^{\delta}_{n(1 - 2 \alpha)} \vert )/n
     \leq \alpha.$ 
Therefore, for all sufficiently large $n,$
\begin{eqnarray}
   \label{lndiffbound}
    \frac{1}{n} \, (\ln \vert \Omega^{\delta}_n \vert - 
   \ln \vert \Omega^{\delta}_{n(1 - 2 \alpha)} \vert ) 
    & = & \frac{1}{n} \, (\ln \vert \Omega^{\delta}_n \vert - 
         \ln \vert \Omega^{\delta'}_{n(1 - 2 \alpha)} \vert ) 
       \nonumber  \\*[0.07cm] & & \, \, + \,
       \frac{1}{n} \, (\ln \vert \Omega^{\delta'}_{n(1 - 2 \alpha)}
       \vert - \ln \vert \Omega^{\delta}_{n(1 - 2 \alpha)} \vert )
      \nonumber  \\*[0.16cm] 
    & \leq & 2 \alpha \ln 4 \, + \, \alpha = \alpha C .
\end{eqnarray}

Now, the bound (\ref{jnlessthanbh}) is larger than the upper bound
for $ \bE_0 ( e^{ - \sip J_n} \, 1_{\{ J_n \leq n^{1- \delta} \}} )$  
in (\ref{thresholdup}) if we choose
\begin{eqnarray}
  \label{condforbstar}
  \ln \vert \Omega^{\delta}_{n(1 - 2 \alpha)} \vert/n + \xi
   - \ln 4 - \sip b & > & 
    \ln \vert \Omega^{\delta}_n \vert /n - \ln 4, 
           \nonumber \\
      \mbox{ equivalently, } \qquad \qquad \qquad 
      \qquad \qquad \qquad \qquad \qquad & &  \nonumber \\ 
        ((\ln \vert \Omega^{\delta}_{n(1 - 2 \alpha)} \vert -
          \ln \vert \Omega^{\delta}_n \vert) /n + \xi )/ \sip & > & b .
\end{eqnarray} 
Note that, since $\xi \geq - \alpha \ln \alpha,$ 
choosing  $\alpha < e^{-C}$ yields  $\xi > \alpha C,$
equivalently, $-\alpha C + \xi > 0.$
Hence, by virtue of (\ref{lndiffbound}), if we choose $\alpha < e^{-C},$ 
it follows that there is 
some $\zeta_* = \zeta_*(\sip) >0$ 
such that for all sufficiently large $n,$
\begin{equation}
   \label{zetachoice}
(\ln \vert \Omega^{\delta}_{n(1 - 2 \alpha)} \vert -
          \ln \vert \Omega^{\delta}_n \vert) /n + \xi > \zeta_*.
\end{equation}  
Consequently, we can let 
\begin{equation}
  \label{littlebstar}
  b_* = b_*(\sip) = \zeta_* / \sip  >0,
\end{equation}    
and, the announced result follows for all $ b < b_*.$ 
Observe that it follows from display (\ref{zetachoice}) that 
$\zeta_*(\sip)$ is decreasing to zero as $\sip \rightarrow 
\infty,$ equivalently, as both $b$ and $\alpha$ tend to zero.
This completes the proof.
\end{proof}

In fact, minor adaptations of the arguments provide a proof
for the case when the bound $n^{1 - \delta}$
in the statement of Proposition \ref{lowerboundforjn} 
is replaced by $n \, q_n,$ where $q_n \rightarrow 0$ arbitrarily
slowly as $n \rightarrow \infty.$  In other words, the set of paths
with $J_n \in [0, n q_n)$ contributes to $ \bE_0 ( e^{ - \sip J_n})$
or to the $k$-th moments $\bE_{\sip}(\chi_n^k)$ merely negligibly
in the sense that the contribution is $o(\bE_0 ( e^{ - \sip J_n}))$
or $o(\bE_{\sip}(\chi_n^k)),$ respectively, as $n$ tends large
(in fact, this error term is exponentially smaller,
as the proof of Proposition \ref{lowerboundforjn} indicates). 
As a consequence of Propositions \ref{upperboundforjn} and
\ref{lowerboundforjn}, for all that follows, we may neglect 
to keep track of those error terms and assume that 
\begin{equation}
  \label{rangejn}
   J_n \in [ b_1 n,  b_2 n ]
\end{equation}
for some constants $0 < b_1 < b_2 < \infty$ 
such that $\sip b_2$ is a positive number independent of $\sip.$


\section{Palm Distribution of the Point Process of Self-Intersections}
\setcounter{equation}{0}  
 
Palm distributions help answer questions dealing with properties
of a point process, viewed from a {\em typical} point. For example, 
we can make mathematically precise the perhaps heuristically 
clear answers to the questions
 (a) what is the mean number of points of a point process in the 
 plane whose nearest neighbors are all at distance at least $r$ ? 
 and (b) what is the probability
that the point process has a certain property,
given that the point process has a sample point at $x$ ? 

We pause for two observations to illustrate the subtlety and
the importance of the concept of Palm distributions.
First, the Palm distribution represents a {\em conditional
distribution} when applied to, say, simple point processes on
the real line, which historically motivated the study of
Palm distributions 
({\sc Palm} \cite{palm} and {\sc Kallenberg} \cite{kall}, 
Chapter 10, p.\ 83). However, the event to be conditioned upon has
probability zero.
Second, notice that, given a realization of a point
process, the sample point closest to the origin is {\em not}
a typical point of the realization but rather special 
since identified as the point with the property ``being closest
to the origin". Consequently, in analyzing typical points 
regarding some properties of the point process,
conditioning upon the point closest to the origin or upon 
any other specified point will lead to an incorrect answer.

Interestingly, Palm distributions may be utilized for other
typical random geometric objects than points. In the study at hand, 
the relevant typical tools will be the 
{\em lines} which are typical relative to the SILT
 $J_n$ when $J_n \in [ b_1 n,  b_2 n ],$
more precisely, relative to the point process
of self-intersections of the SRW with 
$J_n \in [ b_1 n,  b_2 n ]$ in a cone
that is defined via a line through the origin.
This idea will be taken up in Section 3.1.

After this passage of motivation, let us introduce more  
notation as needed. 
If we let $X_n$ and $Y_n$ denote the first and second coordinate
processes of the SRW, that is, $S_n = (X_n, Y_n)$ for every 
integer $n \geq 0,$ define the distance
\begin{equation}
 \label{distance}
  \chi_n = (X_n^2 + Y_n^2)^{1/2} 
\end{equation}  
or the root of the square displacement
of the walk from the starting point ${\bf 0}.$  
Moreover, let $ \bP_{\chi_n}$ denote the probability distribution
of the distance $\chi_n$ of the SRW. On the range of $x$ where we
can invoke the Local Central Limit theorem, we will approximate
$ d \bP_{\chi_n}(x) / dx $ by the density of the corresponding
Brownian motion, that is, $ 2 (2 \pi n)^{-1/2} \exp \{ - x^2/2n \},$
and on the remaining range of $x$ in $(0,n],$ 
use some bounds on a large deviation estimate 
for $d \bP_{\chi_n}(x),$ which denotes
$ \bP_0 ( x \leq \chi_n < x+ dx)$ 
when $dx$ is arbitrarily small. 

Note that throughout the paper, we shall omit discussion of the 
obvious case $\sip =0.$ We next collect a technical lemma that
relies on a condition and a couple more definitions.
The main players in this condition are
bounded numbers $a_x$ that depend on $x$
and are such that there is some number $\zeta >0$ so
that $\sip a_x \geq \zeta$ for every $0 \leq x \leq n.$ 
Then define
\begin{eqnarray}
             \label{mux}
   \mu_x & = & (\sip a_x)^{1/2} \, n^{3/4}  \\
             \label{functionk}
    q(x) & = &  \exp \{ - \sip \frac{a_x}{2}  n^{1/2} \} 
\end{eqnarray}   
for every $n \geq 0,$ $\sip > 0,$ and $x $ in $[0,n].$
Since $a_x $ is bounded in $x,$
for suitably small $\varepsilon \geq 0$ and for $ \gamma>0,$ 
we may define 
\begin{eqnarray} 
  \label{regioninter}
      r_1 & = & r_1(\varepsilon, \gamma) = 
           \sup \{ x \in [0, n] : 
          x \leq \gamma \mu_x n^{- \varepsilon} \} \nonumber \\
      r_2 & = & r_2( \gamma)  = 
           \sup \{ x \in [0, n] :   x \leq \gamma \mu_x  \} . 
\end{eqnarray} 
Thus, $r_2( \gamma) = r_1(0, \gamma).$


\smallskip
{\sc {\bf Condition D.} }
For any suitably small $\varepsilon \geq 0,$ 
there exist some $\gamma >0$ and $ \rho_* > 0$ such that 
\begin{equation}
  \label{conditiond}
   \int_{r_1}^{n} \, x \, q(x) \,  d \bP_{\chi_n}(x)
    = \rho_n \, \int_0^{r_1} \, x \, q(x) \,  d \bP_{\chi_n}(x) 
\end{equation}
with $\rho_n \geq \rho_* $ for all sufficiently large $n.$

\medskip
If $ \int_{0}^{r_2} \, x \, q(x) \,  d \bP_{\chi_n}(x)
= o( \int_{r_2}^{n} \, x \, q(x) \,  d \bP_{\chi_n}(x))$
as $n \rightarrow \infty,$ then $\varepsilon =0$ and $\rho_n
\rightarrow \infty.$
Observe that, by virtue of the expression in (\ref{functionk}) 
for $q(x),$ Condition D guarantees 
that $a_x$ not be constant in $x$ and $\sip >0.$

\begin{lemma}{\bf (Exponent of Expected Radial Distance equals $3/4$)}
    \label{twointegrals}
Let $\sip > 0.$ Assume that the $a_x$ are bounded 
numbers that depend on $x,$ are such that there 
is some number $\zeta >0$ so that $\sip a_x \geq \zeta$ 
for every $0 \leq x \leq n,$ and that satisfy Condition D 
in $(\ref{conditiond})$ for some $\varepsilon \geq 0$ and $\gamma >0.$
Define
\begin{eqnarray}
    \label{threeintegrals}
 I_n  & = &  \int_0^n \, x \, q(x) \,  d \bP_{\chi_n}(x)  \\
 g(n) & = &  \int_0^n \, (a_x)^{1/2} q(x) \,  d \bP_{\chi_n}(x),
                              \nonumber 
\end{eqnarray}
where $q(x)$ is defined in $(\ref{functionk}).$
Then there are some constants $ M < \infty$ and $c(\rho_*) > 0$ 
(both independent of $\sip$) such that as $n \rightarrow \infty,$ 
\begin{equation}
  \label{inint}
  \gamma  \, c(\rho_*) \, \sip^{1/2} \, 
   n^{3/4 - \varepsilon} \, (1+o(1))
  \leq \frac{I_n}{g(n)} \leq  M \, \sip^{1/2} \, 
   n^{3/4} \, (1+o(1)). 
\end{equation}
\end{lemma}

\begin{proof}
Luckily, the crudest of all estimates will serve us.
Condition D will only be relevant to the lower bound for $I_n.$
We start with a number of general observations. 
Importantly, note that we shall not need a concrete form of
the expression for $\bP_{\chi_n}(\cdot)$ to prove the statement
of the lemma (and ultimately, of the main result
on the MSDE). However, in order to bound $I_n,$ it is crucial
to recognize the exact nature of the 
exponential function that appears in the integrand. For this
purpose, we are interested in an estimate for $d \bP_{\chi_n}(x),$  
more precisely, 
in an estimate for $d \bP_{\chi_n}(x)$ 
beyond the range $(0, n^{2/3})$ of $x$ where 
the Local Central Limit theorem is in force.
A large deviation type result of 
{\sc Billingsley} \cite{bill}, Theorem 9.4, p.\ 149
says that if $\tilde{S}_n$ denotes the partial sum of $n$
independent and identically distributed random variables
with mean $0$ and variance $1$ and 
$\alpha_n \uparrow \infty$ denotes any sequence so
that $\alpha_n / n^{1/2} \rightarrow 0$ as $n \rightarrow \infty,$
then $\bP ( \tilde{S}_n \geq \alpha_n n^{1/2}) = 
2 \exp \{ - \alpha_n^2(1 + \xi_n)/2 \}$
for some sequence $\xi_n \rightarrow 0.$ 
Now, consider the {\em diagonal} symmetric simple random walk,
each coordinate of which independently takes values $+1$ and
$-1$ with probability  $1/2,$ convolute the two coordinates,
and, apply {\sc Billingsley}'s estimate to each coordinate
separately. The distance of the diagonal random walk (up
to scaling by $\sqrt{2}$) indicates the distance of the 
(non-diagonal) random walk in the square lattice since
we can turn the lattice by the angle $\pi/4.$ 
Applying these steps, we obtain, for $n^{1/2} < \! \! < x \leq n,$ 
\begin{eqnarray}  
    \label{ldestimatone}  
   \bP_0 ( \chi_n \geq x ) & = &  2  \exp \{ - x^2(1 + \xi_x)/(2n) \}
\end{eqnarray}
with $\xi_x = o(1)$ as $x \uparrow $ (increasing) 
and $n \rightarrow \infty.$
If the distance of the SRW is measured along any {\em fixed} line
through the origin and the endpoint of its path, then the values
of $x$ form a discrete set (for fixed $n$). If this distance
is measured along  {\em all} lines through the origin
and the possible endpoints of the SRW after $n$ steps 
(for fixed $n$), then the values of $x$ form a discrete
set as well. Call it $\cZ_n.$ Between the points in $\cZ_n,$ 
we may interpolate the upper tail probability of $\chi_n$ in any
desired way as long as the estimate in (\ref{ldestimatone}) 
is not violated, thereby introducing an error to $I_n$ which can
be shown to be of order less than $o(g(n))$ 
as $n \rightarrow \infty.$ We think of embedding the set $\cZ_n$ in 
the nonnegative reals and of extending (\ref{ldestimatone}) to
a differentiable function that obeys the expression for the
upper tail probability for some partial sum $\tilde{S}_n$ as
prescribed by {\sc Billingsley}'s estimate. In that case,
for arbitrarily small $dx$ for every $n^{1/2} < \! \! < x \leq n,$
expanding the difference 
$\bP_0 ( x \leq \chi_n < x+ dx)$ into a Taylor series yields
\begin{eqnarray} 
   \bP_0 ( x \leq \chi_n < x+ dx) & = & 
     \bP_0 ( \chi_n \geq x) -  \bP_0 ( \chi_n \geq x+ dx) \nonumber 
           \\[0.05cm]
      & = & 2  \exp \{ - x^2(1 + \xi_x)/(2n) \} \cdot 
                [ 1 - \exp\{ -[ x^2(\xi_{x+dx} - \xi_x)  \nonumber \\
             & & \mbox{} \; \,  \, + 
                (2x (dx)+ (dx)^2)(1+ \xi_{x+dx})] /(2n) \} ] \nonumber
                        \\[0.07cm] 
      & = &   \label{ldestimate}    
         n^{-1/2} \exp \{ - x^2(1 + \xi'_{x,n})/(2n) \} \, dx
\end{eqnarray} 
with $\xi_x = o(1)$ and $\xi'_{x,n} = o(1)$ 
as $x \uparrow $ and $n \rightarrow \infty.$
We let $d \bP_{\chi_n}(x)$ denote the difference
$ \bP_0 ( x \leq \chi_n < x+ dx),$ as specified by (\ref{ldestimate}),
when $dx$ is arbitrarily small. 
Observe that (\ref{ldestimate}) is a valid expression as well when the
Central Limit theorem applies. 

Next, upon completing the square
\begin{eqnarray}
 \frac{x^2}{2n} + \sip \frac{a_x}{2} n^{1/2} 
      & = & \label{complsquare}
          \frac{1}{2n} (x^2 + \mu_x^2 )  
         =  \frac{1}{2n} (x - \mu_x )^2 
              +   \frac{1}{n} \, x \, \mu_x \, 
\end{eqnarray}
and by relying on (\ref{ldestimate}), 
we obtain as $n \rightarrow \infty,$ 
\begin{eqnarray}
     I_n & = &  \label{densrepr}
                  \int_0^n \, n^{-1/2} \, x \, 
                \exp \{ - \frac{x^2}{2n} \} \, 
                      \exp \{ - \sip \frac{a_x}{2} n^{1/2} \} \,
                     \exp \{ - \frac{x^2}{2n} \xi'_{x,n} \} \, dx 
                      \\*[0.2cm] 
                 & = &     
                     \int_0^n \, n^{-1/2} \, x \, 
                 \exp \{ - \frac{(x - \mu_x)^2}{2n} \} \,
                \exp \{ - x \, \frac{\mu_x}{n} \} \,
                  \exp \{ - \frac{x^2}{2n} \xi'_{x,n} \} \, dx .
                  \nonumber 
\end{eqnarray}

{\bf (i) Upper Bound for $I_n.$} 
In view of the form of the expression for the integrand in 
(\ref{densrepr}), the key contribution to the integral 
$I_n$ stems from values of $x$
for which the exponential functions are largest. 
It suffices to regard the first two exponential factors of
the integrand. We shall argue
that these exponential factors are not maximal for $x$ of order
{\em strictly} larger than $n^{3/4}.$

Fix some  $\epsilon >0.$  
The first observation is that 
$\exp \{ - (x - \mu_x)^2/2n \} $
decays rapidly with $x$ for $x \geq n^{3/4 + \epsilon},$
and on that interval, has strictly smaller exponential rate
than for $x \leq n^{3/4 - \epsilon},$ the latter rate essentially 
being $- \mu_x^2/2n.$ Another number of observations will indicate
rapid decay of  
$ e_*(x) = \exp \{ - x \mu_x /n  \} $ 
with $x$ for $x \geq n^{3/4 + \epsilon}.$
Indeed, since $a_x$ is bounded in $x,$
if we let $x \leq n^{3/4 - \epsilon} <  n^{3/4} <  
n^{3/4 + \epsilon} \leq u,$ 
then we collect as $n \rightarrow \infty,$
\begin{eqnarray}
   \label{estarest}
   e_*(u) & < & e_*(n^{3/4}) < e_*(x), \\
   e_*(u) & = & o( e_* (n^{3/4})), \nonumber  \\
   e_*(n^{3/4}) & = & o( e_* (x)). \nonumber
\end{eqnarray} 
For some suitably large $ M < \infty$ 
($M = 1 + 1/\zeta^{1/2}$ should suffice) write
$
T_{\epsilon} =  \inf \{ x \in [0, n] : 
          x > (M-1) \mu_x n^{ \epsilon} \}. 
$
Thus,  $x \in (0, T_{\epsilon})$ implies that
for all sufficiently large $n,$ we have
$x - \mu_x \leq  (M-1) \mu_x  n^{ \epsilon}.$   
Clearly, $T_{\epsilon} \geq \zeta^{1/2} (M-1) n^{3/4+ \epsilon},$
in particular,  $T_{\epsilon} \geq  n^{3/4+ \epsilon}$ if
$\zeta^{1/2} (M-1) \geq 1.$
Putting each of these pieces together and recalling
that the integrand has exponential form yields
as $n \rightarrow \infty,$ 
\begin{eqnarray*}
     I_n & = &  \int_0^n \, 
              (x - \mu_x ) \,   q(x) \, d \bP_{\chi_n}(x) 
                 +   \int_0^n \, 
                    \mu_x \,   q(x) \, d \bP_{\chi_n}(x) 
            \\*[0.15cm]  & = &  
                 \int_0^{T_{\epsilon}} \, 
                  (x - \mu_x ) \,   q(x) \, d \bP_{\chi_n}(x) 
                 \, + \, \int_{T_{\epsilon}}^n \, 
                  (x - \mu_x ) \,   q(x) \, d \bP_{\chi_n}(x) 
                     \, + \, \int_0^n \, 
                    \mu_x \,   q(x) \, d \bP_{\chi_n}(x) 
                  \\*[0.15cm]  & = &  
               (1 +o(1))  
                   \int_0^{T_{\epsilon}} \, 
                  (x - \mu_x ) \,  q(x) \, d \bP_{\chi_n}(x)
                 \, + \, 
                    \int_0^n \, 
                    \mu_x \,    q(x) \, d \bP_{\chi_n}(x) 
                 \\*[0.15cm]  & \leq &   
              (1 +o(1)) \, 
             n^{3/4 + \epsilon} (M-1) \sip^{1/2} \, 
                 \int_0^{T_{\epsilon}} \, 
                        (a_x)^{1/2}   q(x) \,  
                  d \bP_{\chi_n}(x) 
                \\*[0.05cm]
                  & & \mbox{} \; \, + \, 
                  \sip^{1/2} n^{3/4} \, 
                 \int_0^n \, (a_x)^{1/2}  q(x) \,  
                  d \bP_{\chi_n}(x)
                   \\*[0.15cm]  & \leq &   
                 (1 +o(1)) \, 
             n^{3/4 + \epsilon} (M-1) \sip^{1/2} \, 
                 \int_0^n \,  
                 (a_x)^{1/2}  q(x) \,  
                  d \bP_{\chi_n}(x) 
                   \, + \,   
                  \sip^{1/2} n^{3/4} \, g(n)
                  \\*[0.15cm]  & = & 
                   (1 +o(1)) \, 
             n^{3/4 + \epsilon} (M-1) \sip^{1/2} \, g(n)
               \, + \,  \sip^{1/2} n^{3/4} \, g(n). 
\end{eqnarray*} 
Since $ \epsilon >0$ was arbitrary, this is the claimed
{\em upper} bound, as $n \rightarrow \infty,$ 
\begin{eqnarray*}  
  I_n & \leq & M \, \sip^{1/2} \, n^{3/4} \, g(n) \, (1+o(1)). 
\end{eqnarray*}

{\bf (ii) Lower Bound for $I_n.$}
To handle the {\em lower} bound for $I_n,$ we suppose the
instance of Condition D and we split the integrals $I_n$ and $g(n),$
respectively, over the three intervals $[0, r_1], $
$(r_1, r_2)$ and $[r_2,n]$ as follows:
\begin{eqnarray} 
   \label{decomposeintothree}
   I_n & = & \int_0^n \, x \, q(x) \,  d \bP_{\chi_n}(x) =
           J_1(n) + J_2(n) + J_3(n) \\*[0.1cm]
   g(n) & = &  \int_0^n \, (a_x)^{1/2} q(x) \,  d \bP_{\chi_n}(x) =
                  \hat{J}_1(n) + \hat{J}_2(n) + \hat{J}_3(n) .
    \nonumber
\end{eqnarray}
In light of the symmetric roles of $J_2(n)$ and $J_3(n)$ 
in Condition D, we may assume that
$ J_3(n) = o ( J_1(n) + J_2(n))$ as $n \rightarrow \infty$
because otherwise $J_1(n)$ can
be expressed in terms of $J_3(n)$ instead of in terms of $J_2(n)$
and parallel reasoning to the one employed below applies to 
establish the lower bound for $I_n/g(n).$
The exponential form of the integrands implies that,
as $n \rightarrow \infty,$  
$ \hat{J}_3(n) = o ( \hat{J}_1(n) + \hat{J}_2(n)).$
Write $J_2(n) = \rho_n J_1(n) $ (here, we neglect a possible factor 
$(1+o(1))$, compare to (\ref{conditiond})) and 
$\hat{J}_2(n) $ $ = \hat{\rho}_n  \hat{J}_1(n)$ 
for some $\rho_n \geq \hat{\rho}_n > 0.$
In addition, observe that, if $J_1(n) = o(J_2(n))$ 
as $n \rightarrow \infty,$
then we obtain $\rho_n, \hat{\rho}_n \rightarrow \infty$ 
as $n \rightarrow \infty$ (In particular, we can choose 
$c(\rho_*) =1$ below). Thus, this case shall be covered as a special
case in our treatment below. A similar remark is in force in
the already excluded scenario
that $J_1(n) + J_2(n) = o(J_3(n))$ as $n \rightarrow \infty.$
Another consequence of the exponential form of the integrands
is that there exists some $\rho_* > 0$ so that 
$\rho_n , \hat{\rho}_n \geq \rho_*$ for all sufficiently large $n$ 
if and only if 
there exists some $\rho_* > 0$ so that 
$\rho_n  \geq \rho_*$ for all sufficiently large $n.$
Thus, $\rho_n$ is bounded away from $0$ if and only if 
$\hat{\rho}_n$ is bounded away from $0.$
Now, since we assume that Condition D is valid, 
it follows that there is a $c(\rho_*) >0$ such that
$( 1 + 1/\rho_n) / ( 1 + 1/\hat{\rho}_n) \geq c(\rho_*)$
for all sufficiently large $n.$
Keeping these in mind, as $n \rightarrow \infty,$ we arrive at            
\begin{eqnarray} 
  I_n &  = &  (1 + o(1) ) \, 
              \int_0^{r_2} \, x \,  q(x) \, d \bP_{\chi_n}(x) 
                     \nonumber \\*[0.1cm]
       & =  &  (1 + o(1) )  \, ( 1 + 1/\rho_n) \, J_2(n) 
                     \nonumber \\*[0.1cm]
       & \geq &  (1 + o(1) )  \, ( 1 + 1/\rho_n) 
              \int_{r_1}^{r_2} (\gamma \mu_x n^{- \varepsilon}) \, 
                      q(x) \, d \bP_{\chi_n}(x) 
                \nonumber \\*[0.1cm]
       & = &  
               (1 + o(1) )  \,  \gamma \, \sip^{1/2} \, 
               n^{3/4 - \varepsilon} \, ( 1 + 1/\rho_n) \, 
              \hat{J}_2(n)  \nonumber \\*[0.1cm]
       & = &   
               (1 + o(1) )  \,  \gamma \, \sip^{1/2} \, 
               n^{3/4 - \varepsilon} \, ( 1 + 1/\rho_n) \, 
          (\hat{J}_1(n) + \hat{J}_2(n) )  \, ( 1 + 1/\hat{\rho}_n)^{-1} \, 
          \nonumber \\*[0.1cm]
             & \geq &  \label{lowboundcondbii}    
               (1 + o(1) )  \,  \gamma \, c(\rho_*)
       \, \sip^{1/2} \, n^{3/4 - \varepsilon} \, g(n),
\end{eqnarray}                 
as desired. This accomplishes the lower bound and proof.                
\end{proof}

We remark that the function $a_x$ will emerge shortly,
in (\ref{axnumbershalf}) below.

\medskip
{\bf 3.1. Point Process of Self-Intersections and Cones.}
Next, we shall transfer the setting
of {\sc Stoyan, Kendall, and Mecke} \cite{skm}, 
Chapter 4, p.\ 99, to SRW language and describe the
particulars of the point process of self-intersections.
Let $\Phi= \Phi_n = \{ x_1, x_2, \ldots \}$ denote the point process
of {\em self-intersection points} of the SRW in ${\bf Z}^2$ when
$J_n \in [ b_1 n,  b_2 n ].$ Thus, $ \vert \Phi \vert
\in [ b_1 n,  b_2 n ].$ We allow the points $x_i$ of $\Phi$
to have multiplicity and count such a point exactly as many
times as there are self-intersections of the SRW at $x_i.$
Observe that $\Phi$ depends on $n,$
$b_1,$ and $b_2,$ thus, on $\sip$ and that the condition 
$J_n \in [ b_1 n,  b_2 n ]$ imposed upon $\Phi$ moves the 
analysis to the large deviation range of the SRW
and to the right setting for the weakly SAW.
This random sequence of points $\Phi$ in ${\bf Z}^2$ 
may also be interpreted as a random measure. 
Note that $\bE_0 \Phi$ is $\sigma$-finite.
Let $N_{\Phi}$ denote the set of all point sequences, generated
by $\Phi,$ ${\cal N}_{\Phi}$ the point process $\sigma$-algebra 
generated by $N_{\Phi},$ and $\varphi \in N_{\Phi}$ 
denote a realization of $\Phi.$ Formally, $\Phi$ is a measurable
mapping from the underlying probability space into 
$(N_{\Phi},{\cal N}_{\Phi})$ that induces a distribution on
$(N_{\Phi},{\cal N}_{\Phi}),$ the distribution ${\bP_{\Phi}}$ 
of the point process $\Phi.$ In light of the $\sigma$-finiteness 
of $\bE_0 \Phi,$ ${\bP_{\Phi}}$ is a probability measure. 
Also, let $\bE_{\Phi}$ denote expectation relative \mbox{to 
$\bP_{\Phi}.$} 

An important intermediate tool will consist in a (weakly self-avoiding)
process related to the SRW $S_n$ which satisfies condition 
(\ref{rangejn}) and whose {\em one}-dimensional distribution 
of the radial component 
we understand well enough to calculate a rather precise expression  
for its expected distance from the origin. In turn, 
this estimate will lead to upper and lower bounds for the expected
distance of the two-dimensional process, and ultimately, for
the expected distance of the weakly self-avoiding walk.
For this purpose, our interest will revolve around
the {\em self-intersections} of the two-dimensional SRW
near (half-)lines, more precisely, within certain {\em cones,} 
positioned at the origin.

A cone will be described by the cone that contains a certain line.
Thus, let us now introduce the {\em test set} $\cV$ of
{\em lines} $L$ that will be useful. Let $\cV$ denote a set of 
half-lines (that we call `lines', for ease) 
that emanate from the origin, spread
around a circle in a way that we will not exactly specify 
at this point but will depend on our (optimal) choice later on (see 
\mbox{Definition \ref{shapeofv}}) and on $n.$ It will turn
out to be efficient to choose the lines in $\cV$ equally spaced
around the unit circle. The choice of $\cV$ will have strong ties 
with the {\em shape} of the set $\Phi.$ Eventually, nothing
else will be retained about $\cV$ than its cardinality 
$\vert \cV \vert.$
While the description of $\cV,$ in particular, of its size
$\vert \cV \vert$ will be precised further in the proof of
Proposition \ref{msdlower},
no more is needed to handle Proposition \ref{distmainterm} 
below, which presents a result
that is valid, regardless of the number $\vert \cV \vert $ of lines
and of the arrangement of lines. 
The proofs of Propositions \ref{msdupper} and \ref{msdlower}, however, 
will address the issue on how to choose the
lines for $\cV,$ in particular, how many are needed to 
allocate the relevant self-intersection points
of the walk to the corresponding cones.
Next, let us turn to the restriction of the process
$\Phi$ to a line $L$ in $\cV,$ more precisely, to all points that 
lie closer to $L$ than to any other line in $\cV.$ For any
$L \in \cV,$ let the ``cone" $\cC_L$ be defined by
\begin{eqnarray}
  \label{cone}
  \cC_L &  =&  \{ x_i \in \Phi : \mbox{dist}(x_i , L) \leq 
                           \mbox{dist}(x_i , L') \mbox{ for all }
                           L \not = L' \in \cV \}
\end{eqnarray}
with the convention that if equality 
 $\mbox{dist}(x_i , L) = \mbox{dist}(x_i , L')$ holds for two lines
 $L$ and $L'$ and a certain number of points $x_i,$ then half of them
 will be assigned to $\cC_L$ and the other half to $\cC_{L'}.$ 
 Note that no point of $\Phi$ belongs to more than one $\cC_L$ 
 and each point to exactly one  $\cC_L.$ Thus,
 $\vert \cC_L \vert $ equals the number of self-intersections
of the planar SRW $S_n$ in a cone at the origin that contains
the line $L.$ 
Moreover, for any constants $0 < a_1 < a_2 < \infty,$
define the random set
\begin{eqnarray}
  \label{halfline}
       \cL_{1/2}= \cL_{1/2}(\Phi)  & = & 
           \{ L \in \cV: \, 2 \vert \cC_L \vert \in 
              [ a_1 n^{1/2}, a_2 n^{1/2} ] \},
\end{eqnarray}
which depends on $a_1, a_2,$ and $\cV.$ We will choose
$a_1$ and $a_2$ such that $a_1 \sip$ and $a_2 \sip $ 
are positive numbers which are independent of $\sip$ and $n.$

\medskip
{\bf 3.2. Distance along Cones with Order $n^{1/2}$ SILT.}
If $h: \bR \times N_{\Phi} \rightarrow \bR_+$ 
denotes a nonnegative measurable real-valued function and 
$\cL_*(\Phi)$ denotes any subset of lines
in $\cV,$ then
since $\bE_0 \Phi$ is $\sigma$-finite, we may disintegrate
relative to the probability measure $\bP_{\Phi},$ 
\begin{equation}
 \label{disintegration}
    \bE_{\Phi}  \left ( \sum_{L \in \cL_*(\Phi)}
                              h(L, \Phi) \right )
               = \int  \sum_{L \in \cL_*(\varphi)}
                              h(L, \varphi)\,  d \bP_{\Phi} (\varphi)
\end{equation} 
(consult also {\sc Kallenberg} \cite{kall}, p.\ 83, and 
{\sc Stoyan, Kendall, and Mecke} \cite{skm}, p.\ 99). 
For a discussion of some examples of Palm distributions of
$ \bP_{\Phi},$ the reader is referred to the Appendix.  

Next, observe that the conditional distribution
$\bP_{\Phi \vert \chi_n }$ of the 
point process $\Phi,$ given $\chi_n,$ is a function of $\chi_n$ 
and depends on condition (\ref{rangejn}), as explained earlier,
so as to produce realizations that satisfy the requirement
$J_n \in [ b_1 n,  b_2 n ].$ 
Apply formula (\ref{disintegration}) with
\begin{equation}
   \label{choosehfct}
   h(L,\Phi) = \frac{ \exp \{ - \sip \vert \cC_L \vert \}}
                    { \vert \cL_{1/2}(\Phi) \vert} ,
\end{equation}
with $\bP_{\Phi \vert \chi_n }(\varphi \vert x)$ in place of 
$\bP_{\Phi} (\varphi),$ and $\cL_* = \cL_{1/2}$ to
define the numbers $a_x= a_x(\cL_{1/2})$ by
\begin{eqnarray}  
   \label{axnumbershalf}
       \exp \{ - \sip a_x  n^{1/2}/2 \} & = &
          \bE_{\Phi \vert \chi_n} ( \, \vert \cL_{1/2}(\Phi) \vert^{-1}  
                          \sum_{L \in \cL_{1/2}(\Phi)} 
                     e^{- \sip \vert \cC_L \vert  } \vert 
                         \chi_n =x )  \\*[0.1cm]
                 & = & 
             \int_{\bZ^2} \, \vert \cL_{1/2}(\varphi) \vert^{-1}  
                \sum_{L \in \cL_{1/2}(\varphi)} 
                     e^{- \sip \vert \cC_L \vert } \, 
                         d \bP_{\Phi \vert \chi_n }(\varphi \vert x)
\nonumber
\end{eqnarray} 
for $0 \leq x \leq n,$ where we set $ \sum_{L \in \cL_{1/2}} = 0$ 
if $\cL_{1/2} = \emptyset.$
Thus, conditioned on the event $\chi_n =x ,$ the number
$a_x n^{1/2}/2$ may be interpreted as  ``typical" SILT relative
to the lines in $ \cL_{1/2},$ equivalently,
$\exp \{ - \sip a_x  n^{1/2}/2 \}$ represents a
``typical'' penalizing factor with respect to $ \cL_{1/2},$
provided that $\chi_n =x.$
Taking expectation, we arrive at the
expected ``typical'' penalizing factor 
\begin{equation}  
  \label{averageline}
   \bE_0 ( e^{- \sip \philh} ) 
            =   \bE_0 ( \exp \{ - \sip a_{\chi_n}  n^{1/2}/2 \}).
\end{equation}
In the same fashion, we calculate  
\begin{equation}  
  \label{distaverpenweight}
   \bE_0 ( \chi_n \, e^{- \sip \philh} ) 
        = \bE_0 ( \chi_n \, 
            \bE_{\Phi \vert \chi_n} ( \, \vert \cL_{1/2}(\Phi) \vert^{-1}  
                          \sum_{L \in \cL_{1/2}(\Phi)} 
                     e^{- \sip \vert \cC_L \vert  } \vert 
                         \chi_n =x ) ) .
\end{equation}
The proofs of Propositions \ref{msdupper} and \ref{msdlower}
below (see also Definition \ref{coneprocess}) will throw
light on the issue of this particular choice of penalizing
weight. Observe in (\ref{axnumbershalf}), though, that
asymptotically with $n,$ the sum is preserved if lines
were included that have larger SILT than $n^{1/2},$ 
in particular, lines that are typical to the SRW (as
opposed to the weakly SAW) and tend to carry much larger 
SILT. Whence, $\bE_0$ might as well be employed as
the expectation relative to the SRW-paths that
are typical to the weakly SAW, which justifies
(\ref{averageline}) and (\ref{distaverpenweight}).

Our first result collects an expression for
the expected distance $\bE_0 ( \chi_n \, e^{- \sip \philh} ) $
in terms of $g(n)$ as defined in Lemma \ref{twointegrals}.
A parallel derivation will provide an expression for
$\bE_0 ( e^{- \sip \philh} ).$ Ultimately, we will be interested in
the quotient of the two expectations.
To justify the eventual transfer 
of the principal results to the SAW, we shall continue to be 
careful about whether constants in $n$ and/or 
$x$ depend on $\sip$ or not and often indicate this.

\begin{proposition}{\bf (Expected Distance Along Cones with Order 
$n^{1/2}$ SILT)}  
   \label{distmainterm} 
Let $\sip > 0.$ If the $a_x,$ specified in $(\ref{axnumbershalf}),$ 
satisfy Condition D in $(\ref{conditiond})$
for $\varepsilon = 0$ and $\gamma >0,$
then there are some constants
$ 0 < \gamma_* \leq  M < \infty$ (independent of $\sip$ as
$\sip \rightarrow \infty$ and $M$ independent of $\sip >0$ as well) 
such that as $n \rightarrow \infty,$   
\begin{eqnarray*}
  \bE_0 ( \chi_n \, e^{- \sip \philh})
          & = & K(n) \, n^{3/4} \sip^{1/2} g(n)(1+o(1))
\end{eqnarray*} 
for $ \gamma_* \leq K(n) \leq M,$
where $g(n)$ was defined in $(\ref{threeintegrals}).$
\end{proposition}

\begin{proof}
A combination of the observations preceding 
(\ref{distaverpenweight}) together with 
(\ref{axnumbershalf}) and (\ref{averageline}) 
and Lemma \ref{twointegrals}
provides, as $n \rightarrow \infty,$   
\begin{eqnarray}
  \bE_0 ( \chi_n \, e^{- \sip \philh } )
           & = &  \bE_0 (\chi_n \,  \bE_{\Phi \vert \chi_n} ( \,
                       \vert \cL_{1/2}(\Phi) \vert^{-1}  
                          \sum_{L \in \cL_{1/2}(\Phi)} 
                     e^{- \sip \vert \cC_L \vert  } \vert 
                         \chi_n =x )) \nonumber \\*[0.2cm]
               & = &  \label{cesaroaverage}
                   \int_0^n \, x \, \bE_{\Phi \vert \chi_n} ( \,
                       \vert \cL_{1/2}(\Phi) \vert^{-1}  
                          \sum_{L \in \cL_{1/2}(\Phi)} 
                     e^{- \sip \vert \cC_L \vert  } \vert 
                         \chi_n =x ) \, d \bP_{\chi_n}(x)
                         \nonumber \\*[0.2cm]
               & = &  
                      \int_0^n \, x  \, (  
                     \, 
                      \int_{\bZ^2} \, 
                      \vert \cL_{1/2}(\varphi) \vert^{-1}  
                        \sum_{L \in \cL_{1/2}(\varphi)} 
                  e^{- \sip \vert \cC_L \vert  } \, 
                         d \bP_{\Phi \vert \chi_n}(\varphi \vert x) )
                           \, d \bP_{\chi_n}(x)
                        \nonumber    \\*[0.2cm]
               & = &  \label{evalpalmone}  
                   \int_0^n \, x  \, 
                            \exp \{ - \sip a_x n^{1/2}/2 \}
                         \, d \bP_{\chi_n}(x)
                          \\*[0.2cm]
                  & = &  \label{gnintegral} 
                   \int_0^n \, x  \, q(x) \, d \bP_{\chi_n}(x)
                          \\*[0.2cm]
                   & = &  K(n) n^{3/4} \sip^{1/2} g(n)(1+o(1))  
\end{eqnarray} 
for $ \gamma_* \leq K(n) \leq M,$ 
where to obtain the last two lines of the display, we apply 
\mbox{Lemma \ref{twointegrals},} 
with $\gamma_* = \gamma c(\rho_*), $ $\varepsilon=0,$ and with
the $a_x$ being bounded and such that there 
is some number $\zeta >0$ so that $\sip a_x \geq \zeta$ 
for every $0 \leq x \leq n.$ These two properties of $a_x$
may be seen as follows. First, since, by (\ref{halfline}),
$ 2 \vert \cC_L \vert / n^{1/2}$ is in $  [a_1, a_2 ],$
the average of the exponential 
terms $ \exp \{- \sip \vert \cC_L \vert \}$ over all lines
in $\cL_{1/2}(\varphi)$ may be rewritten as 
$\exp \{ - \sip a_x $ $ n^{1/2}/2 \},$ say,
for some number $a_x \in  [a_1, a_2 ],$ 
depending on $x.$ 
In particular, the $a_x$ are bounded. 
Additionally, we assumed (remark following (\ref{halfline}))
that $a_1 \sip$ is a positive number independent of $\sip,$
thus, there is some number $\zeta >0$ so that $\sip a_x \geq \zeta$
for all $x.$ 
This completes our proof.
\end{proof}

\medskip
{\bf 3.3. Weakly Self-Avoiding Cone Process relative to
$r$-Shaped $\Phi.$}
Once the lines are selected for $\cV,$ 
we may classify them according to the SILT that their cones 
carry. For any suitably small $\delta >0,$ define
\begin{eqnarray}
  \label{fourlines}
   \cL_{1/2 \pm } & = & \cL_{1/2 \pm }(\Phi)  =  
        \{L \in \cV: \, 2 \vert \cC_L \vert
                 \in [ a_1 n^{1/2 - \delta},  a_2 n^{1/2 + \delta}] \} \\
   \cL_{-} & = & \cL_{-}(\Phi)  = 
           \{ L \in \cV: \, 2 \vert \cC_L \vert 
                    \in (0,  a_1 n^{1/2 - \delta}) \} \nonumber \\
    \cL_{+}  & = & \cL_{+}(\Phi)  =  
             \{ L \in \cV: \, 2 \vert \cC_L \vert
                     \in ( a_2 n^{1/2 + \delta}, 2 b_2 n] \}   \nonumber \\ 
   \cL_{r} & = & \cL_{r}(\Phi) = 
              \{ L \in \cV: \, 2 \vert \cC_L \vert
              \in [ a_1 n^{r}, a_2  n^{r} ] \}  
                          \nonumber \\  
   \cL_{\emptyset} & = & \cL_{\emptyset}(\Phi) = 
              \{ L \in \cV: \, \vert \cC_L \vert = 0 \} 
                   \nonumber
\end{eqnarray} 
for each $ 0 \leq r \leq 1$ and the same constants 
$ 0 < a_1 < a_2 < \infty$ as employed in (\ref{halfline}).
Thus, we here modify and
extend the earlier definition $ \cL_{1/2}. $ 

In dealing with the problem to derive the
expected distance with respect to the measure $\bQ_n^{\sip},$
in other words, the expected distance of the weakly SAW,
we will introduce a weakly self-avoiding
process that is related to the weakly SAW.
This related object that we shall construct is suitable to 
calculate concrete expressions for the expected 
distances and attempts to ``mimic" the following idea
to asymptotically calculate the expected 
distance of the SRW after $n$ steps
from the starting point (for which process, though, the 
calculation is much more straightforward). In case of the
latter, the basic ingredients may be sketched as follows.
Rely on the Local Central Limit theorem and 
rewrite the density of the approximating Brownian motion to the
SRW in polar coordinates. Calculating the expected distance
of the SRW involves an integration over the {\em radial part}
and an integration over the {\em angle.} This approximation by
means of Brownian motion involves controlling an error.

In case of the former, roughly speaking,
the process may be depicted as a weakly self-avoiding process whose
penalizing weight takes into consideration the number of
self-intersections near the line that passes through the starting 
point and the endpoint of the SRW-path (rather than penalizing
the two-dimensional process according to $J_n$). Importantly,
the definition of this process will depend on 
the choice of the set $\cV,$ as
made precise shortly, which will determine the SILT near the
relevant lines. Moreover, bounds on the expected distance
of this newly-defined process will be gotten by \\*[0.05cm]
\mbox{} \,
(a) keeping track of the radial part of the SRW, penalized by
the SILT in a certain cone,
\\*[0.05cm] \mbox{} \,
(b) by integrating out over all lines in $\cV.$  \\*[0.05cm]
Part of our strategy involves relating the expected distance
of this process with the one of the weakly SAW.
We begin to describe the ``shape'' of the set $\cV.$
Note that $\cV$ depends on $\Phi$ and its so-called shape
reflects upon the shape of $\Phi.$

\begin{definition}[{\bf $\cV$ and $\Phi$ are $r$-shaped}]
  \label{shapeofv}
Let $\rho >0$ be suitably small.
We say that $\cL_r$ {\em contributes (to $J_n$) essentially}
if
$$ 
    \sum_{L \in \cL_r} \vert \cC_L \vert \geq  \frac{1}{2} \,
                     J_n^{1- \rho}. 
$$
In this case, we say that $\cV$ and $\Phi$ are {\em $r$-shaped} or
have {\em shape $r.$}
In particular, when $r=1/2,$ then we say that $\cV$ and $\Phi$ 
have {\em circular shape} or are {\em circular.}
The convention is that multiple shapes
are allowed, that is, $\Phi$ may simultaneously 
have shape $1/2$ and shape $3/4.$
\end{definition}
{\bf Remarks.} \\
{\bf (1)} For our purposes and later calculations, 
it is not necessary that the lines contributing essentially,
as explained in Definition \ref{shapeofv}, have exact SILT
of order $n^r$ in the sense that the real value $r$ is hit
precisely. Instead, it suffices to replace $\cL_r$ by
$\cL_{r \star} =
 \{L \in \cV: \, 2 \vert \cC_L \vert
      \in [ a_1 n^r,  a_2 n^{r + \delta}] \} $ for $\delta >0,$
and to ultimately let $\delta \rightarrow 0$ in the obtained
results (because $\delta >0$ was arbitrary). Hence, when
applying Definition \ref{shapeofv}, we will think of
$\cL_{r \star} $ rather than $\cL_r$ and refer to
\begin{equation}
  \label{approxshape} 
   \sum_{L \in \cL_s \atop \mbox{for } r \leq s \leq r + \delta } 
    \vert \cC_L \vert \geq  \frac{1}{2} \, J_n^{1- \rho} . 
\end{equation} 
With this meaning, it is obvious that, 
for sufficiently large $n,$ there must be $0 \leq r \leq 1$
such that the set $\cL_{r \star} $ contributes essentially, and thus, 
the shape of $\Phi$ and $\cV$ is well-defined.
Nevertheless, for the sake of not complicating our presentation,
we shall not write $\cL_{r \star} $ and not use the extension
in (\ref{approxshape}) but simply write $\cL_r.$ \\*[0.1cm]
{\bf (2)} 
We might as well choose $J_n \tau_n /2$ with $\tau_n \rightarrow
0$ arbitrarily slowly as $n \rightarrow \infty$ in place of 
$J_n^{1- \rho}/2$ in the defining inequality for the shape of $\Phi.$ 
There is nothing special about the choice above.

\smallskip
Observe that if $\cL_r$ contributes essentially then,
by (\ref{rangejn}) and (\ref{fourlines}), 
\begin{eqnarray}
  \label{boundsonlr}                     
   \frac{b_1}{a_2} \, n^{1-r-\rho} & \leq & 
        \vert \cL_r \vert \leq \frac{2 b_2}{a_1} \, n^{1-r} .
\end{eqnarray}
It is apparent that the upper bound in (\ref{boundsonlr})
holds even when $\Phi$ is not $r$-shaped. Since we choose 
$a_1$ and $a_2$ such that $\sip a_1$ and $\sip a_2$ are
independent of $\sip,$ it follows that $b_2/a_1$ is
independent of $\sip.$ Next, similarly as in (\ref{axnumbershalf}), 
for any subset $ \cL$ of $ \cL_r \subset \cV,$ 
define the numbers $a_x = a_x(\cL)$ by
\begin{eqnarray}  
   \label{axnumbershalfr}
       \exp \{ - \sip a_x (\cL)  n^{r}/2 \} & = &
          \bE_{\Phi \vert \chi_n } ( \, \vert \cL (\Phi) \vert^{-1}  
                          \sum_{L \in \cL(\Phi)} 
                     e^{- \sip \vert \cC_L \vert  } \vert 
                         \chi_n =x )  
\end{eqnarray} 
for $0 \leq x \leq n,$ where we set $ \sum_{L \in \cL} = 0$ 
if $\cL= \emptyset,$ and
in parallel to (\ref{averageline}) and (\ref{distaverpenweight}),
define the expected ``typical'' penalizing factor with respect to 
$\cL \subset \cL_r$ by
\begin{equation}  
  \label{averageliner} 
   \bE_0 ( e^{- \sip J_n^{\cL}})
            =   \bE_0 ( \exp \{ - \sip a_{\chi_n}  n^{r}/2 \})
\end{equation}
and  $\bE_0 ( \chi_n \, \exp \{- \sip  J_n^{\cL} \}). $

\begin{definition}[{Weakly self-avoiding cone process relative 
to $r$-shaped $\cV$}]
  \label{coneprocess}
Define a {\em weakly self-avoiding cone process relative 
to $\cV$ in shape $r$} by some two-dimensional process whose
radial part is induced by the probability measure 
\begin{equation}
  \label{distconeprocess}
    \bQ^{\sip, \cV, r}_n = \frac{\exp \{ - \sip \vert \cC_L \vert \}}{
                       \bE_0 \exp \{ - \sip  J_n^{\cL_r}  \}} 
\end{equation}
on the set of SRW-paths of length $n$
if $\cV$ has shape $r,$ where $L$ denotes the line through 
the origin and the endpoint of the SRW after $n$ steps.
Moreover, the expectation  
 $\bE_{\sip, \cV, \cL_r} = \bE_{\bQ^{\sip, \cV, r}_n}$
relative to the radial part is calculated as in 
$(\ref{axnumbershalfr})$ followed by $(\ref{averageliner})$ 
with $\cL = \cL_r.$
\end{definition}

Let  $\bE_{\sip, \cV, *(r)}$ denote expectation of the 
{\em two}-dimensional weakly self-avoiding cone process relative 
to $\cV$ in shape $r.$ In particular, we write 
$\bE_{\sip, \cV, *} = \bE_{\sip, \cV, *(1/2)}.$ 
Thus, the definition of this process depends on the choice
of $\cV$ and on $\Phi.$
Note that there is no unique such process since only 
the distribution of the radial component of the process is
prescribed and not even the distribution on the lines in $\cV$
is specified. Consequently, there will be several ways to choose the
set $\cV.$  Importantly though, the shape carries much information.
It is worthwhile noting that,
while $\bE_{\sip} (\chi_n) $ does not easily appear 
to be accessible to direct calculations,
rather precise asymptotic expressions may be calculated 
for $\bE_{\sip, \cV, \cL_r}(\chi_n)$ for $0 \leq r \leq 1.$

\begin{lemma}{\bf (The $a_x(\cL_{1/2})$ satisfy Condition D)} 
   \label{conditiondsatisfied}
If $\Phi$ has circular shape for sufficiently large $n,$
then the $a_x(\cL_{1/2}),$ 
defined in $(\ref{axnumbershalf}),$ 
satisfy Condition D in $(\ref{conditiond})$
for $\varepsilon = 0$ and $\gamma >0,$ independent of $\sip$
as $\sip \rightarrow \infty.$
\end{lemma}

\begin{proof}
Fix some suitably small $\varepsilon >0$
and suppose that $\Phi$ be circular for all sufficiently large $n.$
Choose $\rho >0$ sufficiently small.
Let us invoke the notation that we introduced in the proof of
Lemma \ref{twointegrals}, that is, write
$ \bE_0 ( \chi_n \, e^{- \sip \philh } ) = I_n
=J_1(n) + J_2(n) + J_3(n),$ and in the same spirit,
$ \bE_0 ( e^{- \sip \philh } )  
= \tilde{J}_1(n) + \tilde{J}_2(n) + \tilde{J}_3(n).$  
We need to show that there is some $\rho_* >0$ so that
$J_2(n) + J_3(n) = \rho_n J_1(n)$ with $\rho_n \geq \rho_* $ 
for all sufficiently large $n>0.$
We begin with proving that $J_2(n) + J_3(n) \not = o(J_1(n))$
as $n \rightarrow \infty.$ 

For a moment, let us suppose in contrast
that $J_2(n) + J_3(n)  = o(J_1(n))$ as $n \rightarrow \infty$
so as to take this claim to a contradiction.
Thus,  $J_2(n)  = o(J_1(n))$ and $J_3(n)  = o(J_1(n))$ 
as $n \rightarrow \infty.$ It would
follow that $I_n = J_1(n)(1 + o(1))$ as $n \rightarrow \infty$
as well as 
$\sum_{i=1}^3 \tilde{J_i}(n) = \tilde{J}_1(n) (1 + o(1)).$   
The probability measure $ \bQ^{\sip, \cV, 1/2}_n $
induces a one-dimensional process $W_n$ which has expectation
$\bE_{\sip, \cV, \cL_{1/2}}( \chi_n),$ call it 
$ \bE_{\sip, \cV, \cL_{1/2}}^{W}(\chi_n). $ 
Associate $W_n$ with the numbers $a_x(\cL_{1/2}).$ 

In view of the exponential form of the integrand of $I_n,$ 
our assumption would imply that there is a number
$z_n=z$ in $  [0, r_1]$ that enjoys the property 
\begin{equation}
   \label{saddleone}
   \frac{ \bE_0 ( \chi_n \, e^{- \sip \philh }) }
        {\bE_0 ( e^{- \sip \philh } ) }  = 
    \frac{I_n}{ \bE_0 ( e^{- \sip \philh } ) }  = 
    (1 + o(1)) \, z 
\end{equation}
as $n \rightarrow \infty.$ In that event,
the function $a_x$ is minimal at $z=z_n,$ 
that is $a_z = \inf_{0 \leq y \leq r_1} a_y$
for all sufficiently large $n.$
This can be seen as follows. 
Define $k(x)   =  \exp \{ - (x^2 + \mu_x^2 )/(2n) \},$   
let $a_0 >0$ and let $0 <\omega \leq a_0$ 
be some arbitrarily small number. If $a_{x_1} = a_0$ and
$a_{x_2} = a_0 - \omega \geq 0$ for $0 \leq x_1, x_2 \leq r_1,$ 
then it follows that $k(x_1) < k(x_2)$ for all
sufficiently large $n.$

Now, for some suitably small $\omega = \omega(\sip) >0,$ define the set 
$$ 
  S_{\omega} = \{ x \in [0, r_1] : a_x > a_z + \omega \} .
$$ 
Consider a modified process $\tilde{W}_n$ that is associated
with numbers $\tilde{a}_x$ with
$\tilde{a}_x = a_x$ for $x \in [0, r_1] \setminus S_{\omega},$
$\tilde{a}_x = a_z + \omega$ for $x \in S_{\omega},$ and
$\tilde{a}_x = a_x + a(n)$ for $r_1 <  x \leq n,$ where
$a(n)>0 $ is some suitable number, chosen so as to preserve the
distribution of $J_n.$ 
Thus, $\tilde{a}_x \leq a_z + \omega$ for $x \in [0, r_1].$  
Observe that the modified process $\tilde{W}_n$ has the same
expectation $ \bE_{\sip, \cV, \cL_{1/2}}^{\tilde{W}}(\chi_n) =
\bE_{\sip, \cV, \cL_{1/2}}^{W}(\chi_n)$
as the process $W_n$ since, firstly, $q(x)$ in (\ref{functionk}) 
was decreased on $(r_1, n],$
and thus, $J_2(n) + J_3(n)$ and the corresponding part
$\tilde{J}_2(n) + \tilde{J}_3(n)$ of the integral
in the denominator of $\bE_{\sip, \cV, \cL_{1/2}}^{W}(\chi_n)$
were both decreased, and secondly, $J_1(n)$ is as before thanks to 
(\ref{saddleone}). Note that adding 
a constant number of self-intersections to {\em all} 
realizations of this underlying weakly self-avoiding process 
$\tilde{W}_n$ does not change its probability measure.
Subtract the number $a_z$ from $\tilde{a}_x$ 
for every $0 \leq x \leq n,$ that is,
let $\hat{a}_x = \tilde{a}_x - a_z \geq 0 $ 
for every $0 \leq x \leq n.$
Thus, $\hat{a}_x \leq \omega$ for $x \in [0, r_1]$
and  $\hat{a}_x$ is suitable on $(r_1, n].$
In particular, we may choose $\omega < a_1 b_1/b_2,$
where $a_1$ was introduced in (\ref{halfline}). 
The gotten process $\hat{W}_n$
associated with the numbers $\hat{a}_x$ has 
expectation  $ \bE_{\sip, \cV, \cL_{1/2}}^{\hat{W}}(\chi_n) =
\bE_{\sip, \cV, \cL_{1/2}}^{W}(\chi_n), $ too, the same
as do $W_n$ and $\tilde{W}_n.$

The number of lines in $\cL_{1/2}$ that would be needed
to assign the self-intersection points of the two-dimensional
process with marginal $\hat{W}_n$ is at least
$ n^{1/2-\rho } b_1 / \omega,$ where $\omega$ is 
suitably small. But since $\rho > 0$ was arbitrary and also by
Remark (2) following Definition \ref{shapeofv}, this contradicts 
(\ref{boundsonlr}) and the assumption that $\Phi$ is circular.
We conclude that $J_2(n) + J_3(n) \not = o(J_1(n))$ as 
$n \rightarrow \infty.$ 
Since $\varepsilon >0$ was arbitrary,
it follows that, for every $\varepsilon >0,$
$J_2(n) + J_3(n) \not = o(J_1(n))$ as 
$n \rightarrow \infty.$

It remains to be shown that there is no subsequence $n_k$
such that $J_2(n_k) + J_3(n_k) = o(J_1(n_k))$
as $k \rightarrow \infty.$ From this it will follow that
there is some number $\rho_* >0$ that bounds $\rho_n$
from below with $n.$ But the same point can be made as
explained above when $n$ is replaced by $n_k$ everywhere.
Whence, we conclude that 
Condition D must hold for $\varepsilon =0.$
Observe that this implies that $r_1 = r_2$ and $J_2(n)= 0.$
 
In addition, we remark that $\gamma >0$ may be chosen 
uniformly over $\sip >0 $ as $\sip \rightarrow \infty.$
This can be seen as follows.
Any of the asymptotic statements in Lemma  
\ref{twointegrals} and in the above lines of proof 
depend on expressions, for example, 
of the form $\sip^{1/2} n^{3/4}.$
Hence, if $N(\sip)$ is a threshold so that, for all
$n \geq N(\sip),$ a given expression in $n$ differs
from its corresponding limiting expression
by at most $\varepsilon$ (some fixed $\varepsilon$), it follows
that $N(\sip') \leq N(\sip)$ for $\sip < \sip'.$ As a
consequence of the fact that $\gamma >0$ may be chosen
uniformly in $n,$ the choice of $\gamma$ is uniformly
over $\sip >0$ as $\sip \rightarrow \infty$
(yet not as $\sip \rightarrow 0.$)
This proves the advertized claim.
\end{proof}

The following two propositions collect the principal results.

\begin{proposition}{\bf (Upper Bound for $ \bE_{\sip} \chi_n$)} 
   \label{msdupper}
Let $\sip >0.$ There is some constant $ M < \infty$ 
(independent of $\sip$) so that as $n \rightarrow \infty,$  
\begin{eqnarray*} 
  \bE_{\sip} (\chi_n) & \leq & 
         (1 + o(1)) \, 
           \frac{   \max_{\cL \subset \cL_{1/2}}  \,
             \bE_0 ( \chi_n \, e^{- \sip J_n^{\cL}} )}
           {\bE_0 ( e^{- \sip  J_n^{\cL_{1/2}}}  )} 
                \\*[0.1cm]
       & \leq & M \, n^{3/4} \, \sip^{1/2} \,  
                   (1+o(1)) \, 
              \frac{  \max_{\cL \subset \cL_{1/2}}  \, g(n)}{h(n)}       
                   \nonumber \\*[0.1cm] 
              & = &  M \, n^{3/4} \, \sip^{1/2} \,  
                    (1+o(1)) \, 
                    \frac{  \max_{\cL \subset \cL_{1/2}} \,
                    \int_0^n \, (a_x(\cL))^{1/2}  q(x) \, 
                    d \bP_{\chi_n}(x)}
                   { \int_0^n \, q(x) \, d \bP_{\chi_n}(x)},
\end{eqnarray*}
where $q(x)$ is as defined in $(\ref{functionk})$ and $a_x(\cL)$ 
in $(\ref{axnumbershalfr})$ when $r=1/2.$
\end{proposition}
 
\begin{proof}
It will be sufficient to prove that, for $\cV$ in circular shape,
as $n \rightarrow \infty,$
\\*[0.1cm] \mbox{} \qquad 
(I) \, \, $ \bE_{\sip, \cV, *} (\chi_n) 
            \leq  (1 + o(1)) \,   
            \max_{\cL \subset \cL_{1/2}}  \,
             \bE_0 ( \chi_n \, e^{- \sip J_n^{\cL}} ) /
            \bE_0 ( e^{- \sip  J_n^{\cL_{1/2}}}  ) , $
\\*[0.05cm] \mbox{} \qquad 
(II) \,  $  \bE_{\sip}(\chi_n)  \leq  \bE_{\sip, \cV, *}(\chi_n)
              (1+o(1)),$ and, 
\\*[0.05cm] \mbox{} \qquad 
(III) to evaluate the expression on the righthand side in (I).

\smallskip
{\bf Parts (I) and (III).}
For this purpose, assume that $\cV$ is $1/2$-shaped
for all sufficiently large $n.$ 
Fix some suitably small $\delta >0$ and fix $\rho < \delta/2.$
Since upon the assignment of all self-intersection points to cones,
every line in $\cV$ falls in exactly one of the
four sets  $\cL_{1/2 \pm},$ $\cL_{-},$  $\cL_{+},$ and 
$\cL_{\emptyset},$ defined in (\ref{fourlines}) and
in view of the definition of $\bE_{\sip, \cV, *},$ we collect
\begin{eqnarray}
    \label{threepalms}
 \bE_{\sip, \cV, *}(\chi_n) 
               & \leq &   
                \sum_{\tilde{\cL} \in \{ 
                  \cL_{1/2 \pm}, \cL_{-}, \cL_{+}, \cL_{\emptyset} \}} 
                 \;  \bP_{\Phi} ( L \in \tilde{\cL})
              \,  \frac{   \max_{\cL \subset \tilde{\cL} }  \,
             \bE_0 ( \chi_n \, e^{- \sip J_n^{\cL}} )}
           {\bE_0 ( e^{- \sip  J_n^{\tilde{\cL}}}  )} ,
\end{eqnarray}
where $ \bE_0 ( \chi_n \, e^{- \sip J_n^{\cL}} )$ and
$\bE_0 ( e^{- \sip  J_n^{\tilde{\cL}}}  )$ 
for $ \tilde{\cL} \in \{\cL_{1/2 \pm}, 
\cL_{-}, \cL_{+}, \cL_{\emptyset}\}$
are to be understood in the sense of definition (\ref{averageliner}).
Observe that the last term in (\ref{threepalms}),
the SRW term (because it resembles the contribution that we would
obtain from the SRW), is of asymptotic order no larger than $n^{1/2}.$  
It will turn out that (\ref{threepalms}) is bounded above by
the first term in (\ref{threepalms}) times $(1+o(1))$ as 
$n \rightarrow \infty.$ Clearly,  
$ \bP_{\Phi} ( L \in \cL_{1/2 \pm}) \leq 1.$ Thus, 
the first term is bounded above by its quotient.
First, we will see that the first term dominates 
the second and third terms and is of order no larger than
$n^{3/4}.$ 

To analyze $ \bE_{\sip, \cV, \cL_{1/2+r}} (\chi_n) =
  \bE_0 ( \chi_n \, e^{- \sip \, \philpr} ) /
        \bE_0 ( e^{- \sip \, \philpr} ), $
we shall proceed much as we did
to verify Proposition \ref{distmainterm}. 
Let us point out the modifications required in the proofs
of Lemma \ref{twointegrals} and Proposition \ref{distmainterm}. 
For ease of exposition, we will shorter 
write $a_x = a_x(\cL_{1/2+r}).$
In parallel to the handling of 
$I_n$ in Lemma \ref{twointegrals}, we let 
$ q_r(x) =  \exp \{ - \sip a_x n^{1/2+r}/2 \}$ for any
$r \in [-1/2, 1/2]$ and evaluate the integral
$$
  I_n(r)= \int_0^n \, x \, q_r(x) \,  d \bP_{\chi_n}(x)
$$
by proceeding along the same reasoning (as in Lemma \ref{twointegrals}) 
with
$$ 
      \mu_x(r) =  (\sip a_x)^{1/2}  \, n^{3/4 + r/2} 
$$
in place of $\mu_x =  (\sip a_x)^{1/2} \, n^{3/4}$ in (\ref{mux}).
Then we arrive at 
\begin{equation}                                        
   \label{rintegral}
     I_n(r) = \sip^{1/2} \, K_r(n) \, n^{3/4 + r/2} g_r(n)(1+ o(1))
\end{equation}
for $  K_r(n) \leq  M < \infty,$ 
where  
$g_r(n) = \int_0^n \, (a_x(\cL_{1/2+r}))^{1/2}  q_r(x) \,  
d \bP_{\chi_n}(x),$
and in view of Proposition \ref{distmainterm} and 
since, by Lemma \ref{conditiondsatisfied},
the $a_x(\cL_{1/2})$ satisfy Condition D in $(\ref{conditiond})$
for $\varepsilon =0$ and some $\gamma >0,$
we have $0 < \gamma_* \leq K_0(n) = K(n) \leq  M < \infty$
($M$ independent of $\sip$ and $\gamma_*$ independent
of $\sip$ as $\sip \rightarrow \infty$). 
Following the steps in the proof of Proposition \ref{distmainterm}
line by line, with $\cL_{1/2}$ replaced by  $\cL_{1/2+r},$ 
and keeping in mind expression (\ref{rintegral}), 
we obtain for each $r \in [-1/2, 1/2],$ as $n \rightarrow \infty,$
\begin{eqnarray}
   \label{expectrlines} 
  \bE_0 ( \chi_n \, e^{- \sip \, \philpr} ) 
            & = &  \sip^{1/2} \, K_r(n) \, n^{3/4 + r/2} g_r(n)(1+ o(1))
\end{eqnarray} 
and the quotient
\begin{eqnarray}
   \label{secondpalm}   
    \bE_{\sip, \cV, \cL_{1/2+r}} (\chi_n) & = &
      \frac {\bE_0 ( \chi_n \, e^{- \sip \, \philpr} ) }
         { \bE_0 ( e^{- \sip \, \philpr} ) } 
             = K_r(n) \, \sip^{1/2} \, n^{3/4 + r/2} \, 
                     \frac{g_r(n)}{h_r(n)} \, (1+ o(1))
                  \nonumber
                 \\*[0.2cm] & = & 
               K_r(n)  \, \sip^{1/2} \, n^{3/4 + r/2} \, \frac{ \int_0^n  
                   \, (a_x )^{1/2}  q_{r}(x) \, d \bP_{\chi_n}(x)}
                   { \int_0^n \, q_{r}(x) \, d \bP_{\chi_n}(x)} \,
                    (1+o(1)) . 
\end{eqnarray}
Hence, in view of the boundedness of the $a_x(\cL_{1/2+r})$ in $n,$
expression (\ref{secondpalm}), for $r \in [- 1/2, -\delta],$      
is maximal for $r= - \delta,$ as $n \rightarrow \infty.$
As a consequence,
\begin{eqnarray}
    \label{secondpalmsec}  
      \frac{\bE_0 ( \chi_n \, e^{- \sip \philm} )}
           {\bE_0 ( e^{- \sip \philm} )} 
            & \leq &   \frac {\bE_0 ( \chi_n \, e^{- \sip \, \philpr} ) }
         { \bE_0 ( e^{- \sip \, \philpr} ) } 
\end{eqnarray} 
with $r= - \delta.$ However, the righthand side of 
(\ref{secondpalmsec}) is strictly
less than $ \bE_{\sip, \cV, \cL_{1/2}} = $ 
$\bE_0 ( \chi_n \, e^{- \sip \philh} ) /$
             $\bE_0 ( e^{- \sip \philh} ). $
Holding on to (\ref{secondpalm}) with $r=0,$ we conclude that   
\begin{eqnarray}
    \label{secondpalmcomp}   
     \frac{\bE_0 ( \chi_n \, e^{- \sip \philm} )}
           {\bE_0 ( e^{- \sip \philm} )} 
         & < & 
           \frac{\bE_0 ( \chi_n \, e^{- \sip \philh} )}
           {\bE_0 ( e^{- \sip \philh} )} 
           \leq M \sip^{1/2} \, n^{3/4} \, 
           \frac{ g(n)}{h(n)} \, (1+o(1)).  
\end{eqnarray} 
Since the exponents of the terms in the expression on the
leftmost side of (\ref{secondpalmcomp}) are strictly less than
the exponent of the leading term in the expression in the middle, 
even more is true,
namely, as $n \rightarrow \infty,$
\begin{eqnarray} 
 \label{secondpalmcompas}
    \frac{\bE_0 ( \chi_n \, e^{- \sip \, \philm} )}
                           {\bE_0 ( e^{- \sip \, \philm} )}  
        & = &  o \left( \frac{\bE_0 ( \chi_n \, e^{- \sip \, \philh} )}
                           {\bE_0 ( e^{- \sip \, \philh} )} 
                       \right). 
\end{eqnarray} 
(\ref{secondpalmcompas}) continues to hold if the numerator of
the expression on the
lefthand side is maximized over subsets of  $\cL_{-}.$ 
We conclude that the second term in (\ref{threepalms}) is dominated
by the rightmost side of (\ref{secondpalmcomp}).     
To accomplish the upper bound for $ \bE_{\sip, \cV, *}(\chi_n),$ 
it remains to be shown that the second factor of the 
first term in (\ref{threepalms}) 
as well dominates the third term.

We will argue that $\bP_{\Phi} ( L \in \cL_{1/2+r})$ for 
$r \in (\delta, 1/2]$ is small relative to 
$\bP_{\Phi} ( L \in \cL_{1/2}).$
The probability  $ \bP_{\Phi} ( L \in \cL_{1/2+r})$ may be interpreted
as a Palm probability (see (\ref{palmdistexamb}) in the
Appendix), that is,
\begin{equation}
   \label{palmrlines}   
     \bP_{\Phi} ( L \in \cL_{1/2+r})  =
           \frac{ \bE_{\Phi} \sum_{ L \in \cV} 1_{\cL_{1/2+r}}(L) }
                 { \vert \cV \vert }
          = \frac{ \bE_{\Phi} \vert \cL_{1/2+r} \vert}
                 { \vert \cV \vert } .
\end{equation} 
There are, however, at most $b_2 n$ self-intersections to
distribute to cones, each of which carries at least $ a_1 n^{1/2 + r}/2$
self-intersections. Thus, $\vert \cL_{1/2+r} \vert
 \leq  (2 b_2/ a_1) n^{1/2 - r}. $ Therefore on the one hand,
 $\bP_{\Phi} ( L \in \cL_{1/2+r})  \leq  (2 b_2/ a_1 \vert \cV \vert) 
\,  n^{1/2 - r}.$  On the other hand, because we assumed $\Phi$ to
be circular, we have 
 $\bP_{\Phi} ( L \in \cL_{1/2})  \geq  (b_1/ a_2 
\vert \cV \vert) \, n^{1/2 - \rho}.$ Hence, a combination of these
two observations together with (\ref{secondpalm}) yields, for every 
$r \in (\delta, 1/2],$
\begin{eqnarray*} 
  \frac{\bE_0 (\chi_n \, e^{- \sip \, \philpr} )}
                           {\bE_0 ( e^{- \sip \, \philpr} )} \, 
                      \bP_{\Phi}( L \in  \cL_{1/2 +r}) 
          &  < &   \frac{\bE_0 ( \chi_n  \, e^{- \sip \, \philh} )}
                           {\bE_0 ( e^{- \sip \, \philh} )}  \,
                   \bP_{\Phi}( L \in  \cL_{1/2}) 
\end{eqnarray*} 
because the order of the term on the left is at most 
$n^{5/4 - r/2}/\vert \cV \vert $ and the one on the righthand
side is at least $n^{5/4-\rho}/\vert \cV \vert, $ the latter
being {\em strictly} larger than the former since we picked
$\rho < \delta/2,$ and thus,                       
\begin{eqnarray} 
 \label{thirdpalmcomp}
    \frac{\bE_0 ( \chi_n \, e^{- \sip \, \philp} )}
                           {\bE_0 ( e^{- \sip \, \philp} )} \, 
                      \bP_{\Phi} ( L \in  \cL_{+}) 
           & <  &  \frac{\bE_0 (\chi_n \, e^{- \sip \, \philh} )}
                           {\bE_0 ( e^{- \sip \, \philh} )} . 
\end{eqnarray} 
Inequality (\ref{thirdpalmcomp}) continues to hold if 
the numerator of the quotient of the lefthand expression is maximized
over subsets of  $\cL_{+}.$ 
Since the exponents of the terms in the expression on the
lefthand side of  (\ref{thirdpalmcomp}) are strictly less than
the exponent of the leading term on the righthand side, again
as $n \rightarrow \infty,$
\begin{eqnarray} 
 \label{thirdpalmcompas}
    \frac{\bE_0 ( \chi_n \, e^{- \sip \, \philp} )}
                           {\bE_0 ( e^{- \sip \, \philp} )} \, 
                      \bP_{\Phi} ( L \in  \cL_{+}) 
           & = &  o \left( \frac{\bE_0 ( \chi_n \, e^{- \sip \, \philh} )}
                           {\bE_0 ( e^{- \sip \, \philh} )}  \right). 
\end{eqnarray} 
We summarize our progress as follows.
Since $\delta >0$ was arbitrary, combining 
(\ref{threepalms}), 
(\ref{secondpalmcompas}), 
and (\ref{thirdpalmcompas}) provides as $n \rightarrow \infty,$
\begin{eqnarray}
   \label{distancerepresent}
 \bE_{\sip, \cV, *} (\chi_n)  
           & \leq & (1+o(1)) \,
      \frac{   \max_{\cL \subset \cL_{1/2}}  \,
             \bE_0 ( \chi_n \, e^{- \sip J_n^{\cL}} )}
           {\bE_0 ( e^{- \sip  J_n^{\cL_{1/2}}}  )} 
                      \nonumber \\*[0.15cm]
          & \leq &     \label{twoterms}  
            M \, \sip^{1/2} \, n^{3/4} \, (1+o(1)) \,  
           \frac{  \max_{\cL \subset \cL_{1/2}}  \, g(n)}{h(n)} 
\end{eqnarray} 
for $ M < \infty$ (independent of $\sip$). 
This completes the verification of (I) along with
the asymptotic evaluation of its righthand side.

{\bf Part (II).}
We turn to showing (II), which will finish our proof.
Recall that $\bE_{\sip, \cV, *(r)}$ denotes 
expectation of the two-dimensional 
weakly self-avoiding cone process relative 
to $\cV$ in shape $r.$ 
In order to compare $\bE_{\sip, \cV, *(r)}(\chi_n)$ and 
$\bE_{\sip}(\chi_n),$ the strategy will be to show that, for fixed
$J_n  \in [ b_1 n,  b_2 n ],$
the number of SRW-paths with $J_n $ whose point process $\Phi$ 
is $r$-shaped is larger than the number of SRW-paths with $J_n $ 
whose point process $\Phi$ is $s$-shaped (but not $r$-shaped) 
for $1/2 \leq r <s.$
We will continue to show that $\cL_s$ for $0 \leq s < 1/2$
plays a negligible role as well. 
Hence, the measure $\bQ_n^{\sip}$ prefers 
circular shape. In other words, most SRW-paths that
satisfy (\ref{rangejn}) arise from a $\Phi$ that is $1/2$-shaped, 
Finally, we shall compare the centers of mass of the
weakly self-avoiding cone process and the weakly SAW.

{\bf (a) $\Phi$ prefers circular shape.}
Fix $J_n$ (and assume that $J_n \in [ b_1 n,  b_2 n ] $).
Partition the interval $[1/2,1]$ into $R$ subintervals, each
of which has equal length, that is, let
$ 1/2 = r_0 < r_1 < r_2 < \ldots < r_R = 1.$
We are interested in comparing
the number of SRW-paths with $J_n$ whose point process $\Phi$ 
is $r_{k-1}$-shaped to the number of SRW-paths with $J_n$ whose 
point process $\Phi$ is $r_k$-shaped (but not $r_{k-1}$-shaped).
For this purpose, we shall give an inductive argument over $k.$
Pick a SRW-path $\gamma$ of length $n$ with $J_n$ whose
point process $\Phi$ has shape $r_k.$ 
We will show that 
\mbox{(i) associated} with $\gamma,$ there is a large set $F_{\gamma}$ of 
SRW-paths whose realizations of $\Phi$ have shape $r_{k-1},$ and
(ii) two sets $F_{\gamma}$ and $F_{\gamma'}$ are disjoint for 
$\gamma \not = \gamma'.$ 
To see this, we {\em cut and paste} the path $\gamma$ as follows.
Let $P_{\gamma}$ denote the smallest parallelepiped that contains
the path $\gamma$ and let $l_{\gamma}$ denote the largest
integer less than or equal to the length of the long side of
$P_{\gamma}.$  Divide $P_{\gamma}$
into sub-parallelepipeds whose sides are parallel to the sides of
$P_{\gamma}$ by partitioning the two long sides of $P_{\gamma}$
into $n_f$ subintervals in the same fashion whose endpoints are
vertices of the integer lattice and by connecting the
two endpoints of the subintervals that are opposite to each
other on the two sides. Shift each of the sub-parallelepipeds
including the SRW-subpaths contained
by a definite amount between $1$ and $K$ ($K$: some constant) 
along one of the two directions
of the shorter sides of $P_{\gamma}$ and reconnect the SRW-subpaths
where they were disconnected. In doing this, the shifts are
chosen such that the new path $\gamma'$ will have shape
$r_{k-1}$ and the total number of connections needed to reconnect
those subpaths equals a number $C_n$ that is constant in $k.$
Observe that such a choice of shifts exists.
When walking through the new path $\tilde{\gamma},$ 
because of the necessary extra steps to reconnect the subpaths,
the last several steps of $\gamma$ will be ignored. 
Note that this latter number of steps is independent of $k.$
Hence, if the pieces
to reconnect are self-avoiding, then $J_n$ is no larger after 
this cut-and-paste procedure than before.
This is always possible for otherwise we shift apart the 
sub-parallelepipeds such that they are sufficiently separated 
from each other.
Now, either we choose the reconnecting pieces such that
$J_n$ is preserved or we ``shift back'' (along the 
direction of the long sides of the parallelepiped) some or all of
the sub-parallelepipeds so that any two parallelepipeds 
overlap sufficiently to preserve $J_n$ and then 
reconnect the SRW-subpaths where they were disconnected. 
Again, we shift in such a fashion that the
total number of connections needed to reconnect
the subpaths equals $C_n.$
The number of these newly
constructed paths in $F_{\gamma}$ grows at least at the 
order that the number of ways does to choose $n_f$ locations 
(to shift) among $l_{\gamma}$ sites, which is a number larger
than $1$ for all large enough $n.$
Hence, the number of SRW-paths with $J_n$ whose 
point process $\Phi$ is $r_{k-1}$-shaped is larger
than the number of SRW-paths with $J_n$ whose 
point process $\Phi$ is $r_k$-shaped. Since
this argument can be made for every $1 \leq k \leq R$
and the number of SRW-paths with $J_n$ whose 
point process $\Phi$ has shape $r_R=1$ is at least $1,$
it follows that the number of SRW-paths with $J_n$ whose
point process $\Phi$ is $r$-shaped is maximal for $r=1/2.$ 

{\bf (b) It suffices to consider shapes $r$ with $r \geq 1/2.$}
Our next point will be to reason that it suffices to
consider only $\cL_r$ with $r \geq 1/2.$ 
Suppose that $0 \leq r < 1/2.$  In view of 
(\ref{secondpalm}), for every $\epsilon >0,$ we obtain 
$\bE_{\sip, \cV, \cL_r}(\chi_n) = o(n^{1/2+r/2 + \epsilon})$
as $n \rightarrow \infty.$ 
Moreover, we would need of order $n^{1-r} > n^{1/2}$ lines 
in $\cL_r$ to allocate all points of $\Phi.$ 
However, again by (\ref{secondpalm}),
the points of $\Phi$ in the cones of at least a positive fraction
of these $n^{1-r}$ lines are expected (under $\bP_{\Phi}$) to
lie at distance of order strictly larger than 
$\bE_{\sip, \cV, \cL_r}(\chi_n).$ Therefore, we may instead 
use lines in $\cV$ along directions that are about ``orthogonal'' 
to the directions of the lines in $\cL_r,$ that is, 
lines that cross the smallest rectangle that contains
the points of $\Phi$ along the long side of the rectangle.
In other words, we may use lines in $\cL_s$ with $s \geq 1/2.$ 
Consequently, it follows that it is sufficient to restrict 
attention to $r$-shaped $\Phi$ for $1/2 \leq r \leq 1$ and to 
use $\cL_r$ with $r \geq 1/2.$ 

The considerations in (a) above also imply that both probability
distributions decay exponentially fast around their centers of mass. 
Combining this observation with the fact that
the shape of $\Phi$ relates the SILT of the weakly SAW to the one of
the weakly self-avoiding cone process provides that the two probability
distributions asymptotically have the same centers of mass (up to error
terms).
Together with these, the upshot of above passages (a) and (b) is that,
in comparing $ \bE_{\sip}(\chi_n)$ to $\bE_{\sip, \cV, *(r)}(\chi_n)$
for $ 0 \leq r \leq 1,$ it is enough to choose $r=1/2$ and to 
study the expected distance of the weakly self-avoiding cone 
process relative to $\Phi$ when in circular shape. 
Hence, in particular, we are led to
$$
   \bE_{\sip}(\chi_n)  \leq  \bE_{\sip, \cV, *}(\chi_n) (1 + o(1))
$$
as $n \rightarrow \infty.$
This accomplishes the proof of (II), and thus, ends the proof.
\end{proof}

\begin{proposition}{\bf (Lower Bound for $ \bE_{\sip} \chi_n$)} 
   \label{msdlower}
Let $\sip >0.$  There is a constant $m >0$ (that may depend on $\sip$)
such that as $n \rightarrow \infty,$  
\begin{eqnarray*}  
         \bE_{\sip} (\chi_n)  & \geq & 
            \bP_{\Phi} ( L \in \cL_{1/2}) \, 
           \frac{   \min_{\cL \subset \cL_{1/2}}  \,
             \bE_0 ( \chi_n \, e^{- \sip J_n^{\cL}} )}
           {\bE_0 ( e^{- \sip  J_n^{\cL_{1/2}}}  )} 
           \\*[0.1cm]
        & \geq & (1 + o(1)) \, m \, n^{3/4} \, \sip^{1/2} \, 
                \frac{  \min_{\cL \subset \cL_{1/2}}  \, g(n)}{h(n)},       
\end{eqnarray*}
where $q(x)$ is as defined in $(\ref{functionk}),$ $a_x(\cL)$ 
in $(\ref{axnumbershalfr})$ when $r=1/2,$ and $h(n)$ as in
\mbox{Proposition $\ref{msdupper},$} and, the minimum
is over subsets $\cL \subset \cL_{1/2}$ that form a subset 
of $\cV$ that is circular for sufficiently large $n.$
\end{proposition}  

\begin{proof}
In parallel as we argued earlier in part (II)(b) of the proof
of Proposition \ref{msdupper} (to prove that shapes $s$ of
$\Phi$ for $0 \leq s < 1/2$ are negligible), we find  
$\bE_{\sip, \cV, \cL_s}(\chi_n) \leq  \bE_{\sip}(\chi_n)$
for every $0 \leq s < 1/2.$ Moreover, integrating out
over the lines in $\cV$ yields
$\bE_{\sip, \cV, *(s)}(\chi_n) \leq \bE_{\sip, \cV, \cL_s}(\chi_n). $
Consequently, we arrive at 
\begin{equation}
  \label{lowercomp}
      \bE_{\sip, \cV, *}(\chi_n) \leq \bE_{\sip}(\chi_n) .
\end{equation}
Fix some $\rho >0.$
Similarly as (\ref{threepalms}) was bounded above
along with (\ref{secondpalmcompas}),
(\ref{thirdpalmcompas}), and our earlier remark about the SRW term
in (\ref{threepalms})
together with (\ref{secondpalm}) with $r=0,$
the lefthand side of display (\ref{lowercomp}) can be seen to
be bounded below by a sum of terms of which we only keep the 
maximal term. Thus,
\begin{eqnarray}
    \label{maxpalm}
 \bE_{\sip, \cV, *}(\chi_n) 
               & \geq &   
                  \bP_{\Phi} ( L \in \cL_{1/2} )
               \,  \frac{   \min_{\cL \subset \cL_{1/2}}  \,
             \bE_0 ( \chi_n \, e^{- \sip J_n^{\cL}} )}
           {\bE_0 ( e^{- \sip  J_n^{\cL_{1/2}}}  )} ,
\end{eqnarray}
where the minimum is over subsets $\cL \subset \cL_{1/2}$ 
that form a subset of $\cV$ that is circular for sufficiently 
large $n.$ Hence, for the rest of the 
proof, we may assume that $\Phi$ has circular shape.
Observe that our preceding reasoning and findings 
(see proof of Proposition \ref{msdupper}) make it clear 
that the three other terms related to the sets $\cL_{-}, \cL_{+},$ and 
$ \cL_{\emptyset}$ are of smaller order. 
This together with (\ref{boundsonlr}) for $r=1/2$ yields 
$ \vert \cL_{1/2} \vert \geq (b_1/a_2) \, n^{1/2- \rho}.$
From the fact that $\Phi$ prefers circular shape (see 
part (II) of the proof of Proposition \ref{msdupper}) 
and the estimates in (\ref{boundsonlr}) for $r=1/2,$
and, from the fact that we can choose $a_1$ such
that $b_2/a_1$ is independent of $\sip,$ we
also conclude that $\cV$ can be (optimally) constructed to have
size $\vert \cV \vert = v_n n^{1/2}$ for 
$0 <  v_n \leq v_2 < \infty$ for all sufficiently large
$n,$ where $v_2$ is independent of $\sip.$ 
Consequently, we end up with 
\begin{equation}
   \label{halfprob}  
     \bP_{\Phi} ( L \in \cL_{1/2})  =
           \frac{ \bE_{\Phi} \sum_{ L \in \cV} 1_{\cL_{1/2}}(L) }
                 { \vert \cV \vert } > m_* \, n^{- \rho}
\end{equation}
for $m_* = b_1/(v_2 a_2) >0 $ and every sufficiently large $n.$ 
Note that since $b_1 /a_2$ depends on $\sip,$ so does $m_*.$
In light of (\ref{secondpalm}) with $r=0,$ 
(\ref{lowercomp}), (\ref{maxpalm}), and (\ref{halfprob}), 
we collect
\begin{eqnarray}
    \label{lbpalmone} 
 \bE_{\sip} (\chi_n) & \geq & m_* \, n^{- \rho} \, 
             \frac{   \min_{\cL \subset \cL_{1/2}}  \,
             \bE_0 ( \chi_n \, e^{- \sip J_n^{\cL}} )}
           {\bE_0 ( e^{- \sip  J_n^{\cL_{1/2}}}  )} 
          \nonumber \\*[0.1cm]
          & \geq & 
             m_* \, \gamma_* \,  (1 + o(1)) \,
             \sip^{1/2} \, n^{3/4 - \rho} \,  
        \frac{  \min_{\cL \subset \cL_{1/2}}  \, g(n)}{h(n)},   
\end{eqnarray}  
where $\gamma_*$ is independent of $\sip$ as $\sip \rightarrow \infty.$
Since $\rho >0$ was arbitrary, this lower bound
is as announced when $m_* \gamma_* = m,$ which finishes our proof.         
\end{proof}

In this paper, we won't address the issue of 
existence of the limit $\lim_{n \rightarrow \infty}
       n^{- 3/4} \, \bE_{\sip} (\chi_n)$
but only consider its $\limsup$ and $\liminf.$ 
Observe that, since $a_1 \leq a_x \leq a_2$ for all $x,$
we collect, 
for any subset $\cL_*$ in $\cV,$
\begin{eqnarray}
    \label{unifboundquot}
    ( a_1)^{1/2} \leq 
       \frac{\min_{ \cL \subset \cL_*} \, g(n)}{h(n)} \leq
        \frac{ \max_{ \cL \subset \cL_*} \, g(n)}{h(n)}
        & \leq & ( a_2) ^{1/2}.  
\end{eqnarray}

\medskip
{\bf 3.4. Distance Exponents of the Self-Avoiding Walk.}
Propositions \ref{msdupper} and \ref{msdlower} in company with a 
few extra arguments will imply the following result, which 
establishes that the results about the distance exponents 
carry over to the SAW.

\begin{corollary}
  \label{liminfsupmsd}
Let $\sip > 0.$
There are some constants $0 < \rho_1 = \rho_1(\sip) \leq  
\rho_2 < \infty$ ($\rho_2:$ independent of $\sip$)
such that 
$$  
   \rho_1 \leq
     \liminf_{n \rightarrow \infty}
         n^{- 3/4} \, \bE_{\sip} (\chi_n) 
        \leq  \limsup_{n \rightarrow \infty}
       n^{- 3/4} \, \bE_{\sip} (\chi_n) 
         \leq \rho_2.
$$ 
In particular, the planar self-avoiding walk has 
distance exponent $3/4$ and its normalized expected distance 
is bounded above by the constant $\rho_2.$
\end{corollary}

\begin{proof}
As a consequence of Propositions \ref{msdupper} and \ref{msdlower}
in combination with (\ref{unifboundquot}), as $n \rightarrow \infty,$
\begin{eqnarray*}
  m \,  (a_1 \sip )^{1/2} \, (1 +o(1)) & \leq &
   n^{-3/4} \, \bE_{\sip} (\chi_n)  \leq M \, (a_2 \sip)^{1/2} \, (1 +o(1)),
\end{eqnarray*}
where $0 < m $ may depend on $\sip,$ even as $\sip \rightarrow \infty,$
and $M < \infty$ is independent of all $\sip >0.$
Since we can choose $a_1$ and $a_2$ such that  
$\sip a_1$ and $\sip a_2$ do not depend on $\sip$ 
(see the remark following (\ref{boundsonlr})),
we let $\rho_1 =  m  (a_1 \sip )^{1/2}$ and
$\rho_2 =  M  (a_2 \sip )^{1/2},$ 
which establishes what we stated for the weakly SAW.

It is left to notice the following about the results
in Propositions \ref{distmainterm}, 
\ref{msdupper} and \ref{msdlower}:
Suppose that, for some fixed $\varepsilon >0$ and
$\sip_0,$ there be an integer
$N(\sip_0, \varepsilon) = N(\sip_0),$ 
such that, for every $n > N(\sip_0),$
a given expression for finite $n$ is within distance
$\varepsilon$ from its corresponding limiting expression
in Propositions \ref{distmainterm}, \ref{msdupper} and 
\ref{msdlower}.
Since, everywhere in the calculated 
expectations (as exponent or multiplicative factor),
$n$  shows up as $n^s \sip^t$
for some powers $s, t>0,$ the threshold $N(\sip_0)$ is valid
for all $\sip \geq \sip_0.$ In other words,
$\sip' > \sip_0$ implies that $ N(\sip') \leq N(\sip_0).$
Therefore, each of the various thresholds $N(\cdot)$ can be chosen
uniformly for all $\sip \geq \sip_0,$ in particular,
uniformly as $\sip \rightarrow \infty.$  
Together with the fact that
$ \bE_{\sip} (\chi_n) / n^{3/4}$ is bounded above by a constant
that is independent of $\sip$ as $\sip \rightarrow \infty$
for every sufficiently large $n,$ this implies that 
$\limsup_{n \rightarrow \infty} \lim_{\sip \rightarrow \infty}
 \bE_{\sip} (\chi_n) / n^{3/4} \leq \rho_2,$ where $\rho_2$
is a constant that is independent of $\sip.$ In other words,
the limsup of the normalized expected distance of the SAW
is bounded above. In particular, this establishes that $3/4$ 
is an upper bound for the distance exponent of the SAW.

To bound the distance exponent of the SAW from below,
it suffices to let $\sip \rightarrow \infty$ in the double limit
in such a way that $\limsup_{n \rightarrow \infty} \vert 
\ln \rho_1(\sip) / \ln n \vert =0$ and to then exchange
the limit as $\sip \rightarrow \infty$ and the liminf as $n \rightarrow
\infty$ of $ \ln  \bE_{\beta} ( \chi_n) / \ln n.$ 
As a consequence, the distance exponent of the SAW
is no less than $3/4.$ In other words, the self-avoiding
walk has the same distance exponent as does the weakly SAW,
as claimed.
\end{proof}

This completes the proof of Theorem \ref{sawdistance}. 
Observe that we did not bound the normalized expected distance
of the SAW from below, whereas we gave an upper bound that is constant.
The next result accomplishes Theorem \ref{rmsd}.

\begin{corollary}
  \label{variance}
Let $\sip > 0.$ There are some constants $0 < \rho_3 = \rho_3(\sip)
\leq \rho_4  < \infty$ ($\rho_4:$ independent of $\sip$)
such that 
$$
   \rho_3 \leq
     \liminf_{n \rightarrow \infty}
       n^{- 3/2} \, \bE_{\sip} (\chi_n^2) 
        \leq  \limsup_{n \rightarrow \infty}
       n^{- 3/2} \, \bE_{\sip} (\chi_n^2) 
         \leq \rho_4 .
$$ 
In particular, the MSD exponent of the planar SAW equals $3/2$
and its normalized MSD is bounded above by $\rho_4.$
\end{corollary}
 
\begin{proof}
We will only address the arguments that show the statements
for the weakly SAW and refer the reader to the proof of
Corollary \ref{liminfsupmsd} for the aspects of transferring
some portion of the results to the SAW.
Since the pattern of proof is as before, we only list 
a short guide. Carry out Propositions \ref{distmainterm},
\ref{msdupper} and \ref{msdlower}
with $\chi_n^2$ in place of $\chi_n.$
In particular, write 
\begin{eqnarray*}
    \bE_0 ( \chi_n^2 \, e^{- \sip \, \philh} )  & = & 
            \int_0^n x^2 \, q(x) \, d \bP_{\chi_n}(x) \\
                     & = & \int_0^n [(x- \mu_x)^2 + 2 x \mu_x - \mu_x^2]
                         \, q(x) \, d \bP_{\chi_n}(x)
\end{eqnarray*}
and proceed along our earlier lines that derived  
$ \bE_0 ( \chi_n \, e^{- \sip  \philh } )$ and  
$ \bE_{\sip} ( \chi_n ).$ Nothing more than minor modifications
are required to wind up with the advertized results.                       
\end{proof}

\begin{theorem} 
  \label{convexhull}
Let $\sip > 0$ and let $R_n$ denote the radius of the convex
hull of the SRW-path $S_0, S_1, \ldots, S_n.$ Then $R_n$ satisfies
all statements in Corollaries $\ref{liminfsupmsd}$
and $\ref{variance}$ with $\chi_n$ replaced by $R_n.$
\end{theorem}         
 
\begin{proof}
Observe that we are interested in the {\em maximal} distance
of $S_0, S_1, \ldots, S_n$ along any line
rather than the distance of the position of
$S_n$ from the starting point. 
The reflection principle gives the upper bound 
$d \bP_{R_n}(x)/dx \leq 2 \, d \bP_{\chi_n}(x)/dx$
whereas the lower bound $d \bP_{R_n}(x)/dx \geq  d \bP_{\chi_n}(x)/dx$ 
is readily apparent. From this, the results are immediately collectable.
\end{proof}              

\medskip
{\bf Remark (Transition $\sip \rightarrow 0$).} 
The transition $\sip \rightarrow 0$ is quite different
from the transition $\sip \rightarrow \infty.$ 
Let us quickly look at what happens to our derived expressions
when $\sip =0.$ In that fictive case (since the results were
proven under the assumption $\sip >0$),
all terms in (\ref{threepalms}) are of asymptotic order $n^{1/2},$
and so is the term in (\ref{maxpalm}). Because this is drastically
different from the case $\sip >0,$ in which case the asymptotic
order of the largest term is $n^{3/4},$ we observe a discontinuity
of the expected distance measures 
and distance exponents of the weakly SAW 
as $\sip \rightarrow 0.$ 
In contrast, the case $\sip \rightarrow \infty$ behaves as any
case for fixed $\sip.$



\appendix

\section{Appendix: Examples of Palm Distributions}
\setcounter{equation}{0}  
\setcounter{section}{4}

The first example is the one alluded to in (a) of the introductory
paragraph of Section 3, with the ``typical'' 
random objects being points.
We will write down the Palm distribution of the random measure
$\bP_{\Phi}$ in either case, when the underlying point 
process is stationary and when non-stationary. 
The second example will study a sum of some
exponential random variables, with the random objects 
being {\em points,} whereas the third example will look 
at some sum of exponential random functionals when the ``typical'' 
random objects are {\em lines.} 
All examples are in ${\bf R}^d$ for $d \geq 1.$ We borrow the
notation introduced in Section 3.1.

{\bf Example 1: Number of points without nearest neighbors
within distance $r.$}
Let  $\Phi = \{ x_1, x_2, \ldots \}$ denote some point process
in ${\bf R}^d$ so that its expectation $\bE \Phi$ is  
$\sigma$-finite.
Let $B_r(z)$ denote the ball of radius $r>0$ centered at the 
point $z$ in ${\bf R}^d$ and $o$ denote the origin in ${\bf R}^d.$  
Define the set
\begin{eqnarray*}
  Y & = & 
     \{ \varphi \in N_{\Phi} : \vert \varphi \cap B_r(o) \vert =1 \}
    \\  & = &  
       \{ \varphi \in N_{\Phi} : \varphi \cap B_r(o) 
           \mbox{ is a singleton} \} 
\end{eqnarray*}                   
in  ${\cal N}_{\Phi},$ let $B$ denote some Borel set in 
 ${\bf R}^d,$ and 
\begin{eqnarray*}
  h(z, \varphi)  & = &  1_{B}(z) 1_Y (\varphi -z ),
\end{eqnarray*}     
where $1_B(\cdot)$ denotes the indicator function of $B.$
We may think of the condition $  1_Y (\varphi -z )$ as removing
all points from a realization $\varphi$ that have any nearest
neighbors at distance less than $r.$ Keeping these in mind,
we might be interested in evaluating the mean number of points
of $\Phi$ in $B$ whose nearest neighbors are all at distance
at least $r.$ Thus, in light of some version of
formula (\ref{disintegration}), 
\begin{equation}
 \label{disintegratpoint}
    \bE_{\Phi}  \left ( \sum_{z \in \Phi }
                              h(z, \Phi) \right )
               = \int  \sum_{z \in \varphi }
                   h(z, \varphi)\,  d \bP_{\Phi} (\varphi),
\end{equation} 
we obtain 
\begin{equation}
 \label{disintexampleone}
    \bE_{\Phi}  \left ( \sum_{z \in \Phi \cap B} 
                         1_Y (\Phi - z ) \right )
               = \int  \sum_{z \in \varphi \cap B}
                           1_Y (\varphi - z )  
                          \,  d \bP_{\Phi} (\varphi) .
\end{equation} 
If we assume that $\Phi$ is a {\em stationary} point process with finite
nonzero intensity $\lambda$ and $\mu_d$ denotes Lebesgue
measure in  ${\bf R}^d,$ then the {\em Palm distribution} $\bP_o$
(at $o$) of  $\bP_{ \Phi}$ is a distribution on
$(N_{\Phi},{\cal N}_{\Phi})$ defined by
\begin{equation}
     \label{palmdistexamone}
    \bP_o (Y) = \int_{ {\bf R}^d}  \sum_{z \in \varphi \cap B} 
             \frac{ 1_Y (\varphi - z )  \, d \bP_{\Phi}(\varphi) }
                  { \lambda \mu_d(B) } .       
\end{equation}
This formula holds for any $Y \in {\cal N}_{\Phi}$ and any
Borel set $B$ of positive volume. Note that, by the stationarity 
of the point process, the definition does not depend on the
choice of $B.$

On the other hand, if the point process is {\em not} stationary,
then the Palm distribution $\bP_o$
of  $\bP_{ \Phi}$ is gotten by normalizing as follows:
\begin{equation}
     \label{palmdistexamb}
    \bP_o (Y,B) = \frac {
               \int_{ {\bf R}^d}  \sum_{z \in \varphi \cap B} 
              1_Y (\varphi - z )  \, d \bP_{\Phi}(\varphi)
               }{ 
                \int_{ {\bf R}^d} \vert \varphi \cap B \vert 
                 \, d \bP_{\Phi}(\varphi) 
                 }      
\end{equation}
whenever this quotient is well defined.
This definition depends on the choice of $B.$ 

\medskip
{\bf Example 2: Average of exponential random functional from a
point's perspective.}  
As in the previous example, let $\Phi = \{ x_1, x_2, \ldots \}$ 
denote some point process in ${\bf R}^d.$
For some real numbers $s_1 < s_2,$ define the set
\begin{eqnarray*}
  Y & = & 
     \{ \varphi \in N_{\Phi} : \vert \varphi \cap B_r(o) \vert =
                [s_1, s_2] \}
\end{eqnarray*}                   
in  ${\cal N}_{\Phi},$ let $B$ denote some Borel set in 
 ${\bf R}^d,$ let $\sip >0$ denote some fixed parameter, and 
define
\begin{eqnarray*}
  h(z, \varphi)  & = &  1_{B}(z) 1_Y (\varphi -z ) 
             \exp \{ - \sip  \vert \varphi \cap B_r(z) \vert \} .
\end{eqnarray*}    
This functional {\em marks} or {\em weighs} each point according 
to the number of nearest neighbors within distance $r,$ where the
{\em penalizing} weight has exponential form. 
The more points cluster, the less they weigh. Nicely isolated points 
have large weights. 
In fact, marking the points of the point process generates a
so-called {\em marked point process} (see
{\sc Stoyan, Kendall, and Mecke} \cite{skm}, p.\ 105).  
The average of $ h$ over points in $B$ may be
interpreted as the weight of points in $B$ when the number
of nearest neighbors within distance $r$ lies in $ [s_1, s_2].$  
Therefore, in view of (\ref{disintegratpoint}), 
\begin{eqnarray*}
 \label{disintexampletwo}
    \bE_{\Phi} (  \! \! \sum_{z \in \Phi \cap B} 
                         1_Y (\Phi - z ) 
                \exp \{ - \sip  \vert \Phi \cap B_r(z) \vert \}  
                 )  
               & = &  \! \int   \! \! \! \sum_{z \in \varphi \cap B}
                           1_Y (\varphi - z )  
                 \exp \{ - \sip  \vert \varphi \cap B_r(z) \vert \}        
                           d \bP_{\Phi} (\varphi) .
\end{eqnarray*} 
 
\medskip
{\bf Example 3: Average of exponential random functional from a
line's perspective.}  
As in the previous two examples, let $\Phi = \{ x_1, x_2, \ldots \}$ 
denote some point process in ${\bf R}^d.$ Let $\cV$ denote
some test set of lines that depends on $\Phi.$ 
From (\ref{cone}), recall the restriction $\cC_L $     
of $\Phi$ to a neighborhood of $L,$
more precisely, those points of a realization of $\Phi,$ closest
to the line $L$ in $\cV,$ as opposed to other lines in $\cV.$
If  $\cL_* \subset \cV,$ $B$ denotes some Borel set in the
set of all lines,
the constant $\sip >0$ denotes some fixed parameter, and 
\begin{eqnarray*}
  h(L, \varphi)  & = &  1_B(L) 1_{\cL_*}(L) 
             \exp \{ - \sip  \vert \cC_L \vert  \},
\end{eqnarray*}    
then we obtain,  by (\ref{disintegration}),
\begin{equation}
 \label{disintexamplethree}
    \bE_{\Phi}  \left ( \sum_{L \in \Phi \cap \cL_* \cap B} 
                   \exp \{ - \sip  \vert \cC_L \vert  \}
                   \right )       
           = \int  \sum_{L \in \varphi \cap \cL_* \cap B}
              \exp \{ - \sip  \vert \cC_L \vert \}             
              \,  d \bP_{\Phi} (\varphi) .
\end{equation} 
In this sum of $ h$ over lines in $\cL_*,$ lines are
penalized according to the number of points of $\Phi$ near them.

In the setting of this paper,
$\Phi$ denotes the point process
of self-intersections of the SRW (conditioned upon
$ \vert \Phi \vert  \in [ b_1 n,  b_2 n ]$) and the
terms in the summation are associated to the lines in some $\cV.$
The calculations on $\Phi$ are carried out under the condition that
the SRW-path ends at distance $x.$ Hence, in that case,
the point process and
its Palm distribution both depend on $x.$ In other words,
the role of $\Phi$ is being played by $\Phi \vert x.$


\end{document}